\documentclass{amsart}

\usepackage{amssymb, amsmath, amsthm}
\usepackage{mathabx}
\usepackage{graphicx}
\usepackage[utf8]{inputenc}
\usepackage{hyperref}
\usepackage{young}
\usepackage{epstopdf}
\usepackage{mathtools}
\usepackage{stmaryrd}
\usepackage{rotating}
\usepackage{tabularx}
\usepackage{dynkin-diagrams}
\usepackage{subcaption} 
\usepackage{longtable}
\usepackage{rotating}
\usepackage[all]{xy}
\usepackage{dsfont}
\usepackage{enumitem}
\usepackage{cleveref}
\usepackage{multirow}
\usepackage{float}
\usepackage{array}
\usepackage{booktabs}

\usepackage{listings}

\definecolor{codegreen}{rgb}{0,0.6,0}
\definecolor{codegray}{rgb}{0.5,0.5,0.5}
\definecolor{codepurple}{rgb}{0.58,0,0.82}
\definecolor{backcolour}{rgb}{0.95,0.95,0.92}
\lstdefinestyle{mystyle}{
    backgroundcolor=\color{backcolour},   
    commentstyle=\color{codegreen},
    keywordstyle=\color{magenta},
    numberstyle=\tiny\color{codegray},
    stringstyle=\color{codepurple},
    basicstyle=\ttfamily\footnotesize,
    frame=leftline/topline/bottomline/lines
    breakatwhitespace=true,         
    breaklines=true,                 
    captionpos=b,                    
    keepspaces=true,                 
    numbers=left,                    
    numbersep=5pt,                  
    showspaces=false,                
    showstringspaces=false,
    showtabs=false,                  
    tabsize=2,
}
\usepackage{witharrows}

\usepackage[margin=3cm, right=4cm]{geometry}

\newtheorem{theorem}{Theorem}
\newtheorem*{thm*}{Theorem}
\newtheorem{lemma}[theorem]{Lemma}
\newtheorem{proposition}[theorem]{Proposition}

\newtheorem{remark}[theorem]{Remark}

\newtheorem{corollary}[theorem]{Corollary}

\theoremstyle{definition}

\newcommand{\Z}{\mathbb{Z}}
\newcommand{\N}{\mathbb{N}}

\newcommand{\weyl}{\mathcal{W}}
\newcommand{\schub}{\mathcal{S}}

\newcolumntype{C}[1]{>{\centering\arraybackslash}p{#1}}

\newcommand{\code}{\mathsf{code}}
\newcommand{\w}[1]{\llbracket #1\rrbracket}
\newcommand{\I}{I}

\newcommand{\rr}[1]{\mathbf{#1}}
\newcommand{\rrt}[2]{\widehat{\mathbf{#1}}_{#2}}

\newcommand{\one}{\mathds{1}}
\newcommand{\notInTheta}{\Theta'}

\allowdisplaybreaks

\title[On the computation of homology of type A real flag manifolds]{On the computation of homology of type A real flag manifolds}

\author[lambert]{Jordan Lambert}
\email{jordanlambert@id.uff.br}
\address{Department of Mathematics, ICEx, Universidade Federal Fluminense, Volta Redonda-RJ, Brazil}

\author[rabelo]{Lonardo Rabelo}
\email{lonardo@ice.ufjf.br}
\address{Department of Mathematics, Federal University of Juiz de Fora, Brazil}

\thanks{This work was supported by CAPES (Coordination for the Improvement of Higher Level Personnel) and FAPERJ (Carlos Chagas Filho Foundation for Supporting Research in the State of
Rio de Janeiro) no. 010.002602/2019}



\begin{document}

\begin{abstract}
    In this paper, we present a closed, computable formula for the cellular homology coefficients of real flag manifolds associated with split real forms of type A. We demonstrate the process using movements within the code diagram for permutations. Additionally, we compute the third and fourth homology groups and provide generators for the free part up to the sixth homology group.
\end{abstract}

\maketitle

\section{Introduction}

This work aims to investigate the underlying topology of real flag manifolds, with a particular focus on contributing to the computation of topological invariants in specific cases, as exemplified by the works of Casian-Kodama \cite{CK2010, CK2013}, Matszangosz \cite{Mat19}, and He \cite{He2017, He2020}. It continues a series of previous studies by the authors, dedicated to developing an algorithm for computing the cellular homology coefficients of real flag manifolds, as outlined in \cite{RSm19}. The setting of minimal flag manifolds of types A, B, and C --corresponding, respectively, to real, odd orthogonal, and isotropic Grassmannians-- proved advantageous for deriving initial results due to their relatively lower complexity (see \cite{Rab16}, \cite{LR22b}). The subsequent phase involved examining partial flag manifolds of type A. In \cite{LR22a}, a new formula was introduced to simplify computational efforts; however, this was insufficient to address all issues related to the determination of specific signs. In the present article, leveraging the combinatorial framework of symmetric groups developed in \cite{LR22c}, we fully resolve this problem by deriving a closed computable formula for the coefficients of any real partial flag manifold of type A (see Theorem \ref{thm:coeffA}).

For any flag manifold, it is well-established that the coefficients are either $0$ or $\pm 2$, necessitating a focus on two principal aspects of this computation. The first is determining when the coefficient is non-zero, and, if so, the second is establishing the correct sign. In \cite{LR22a}, the first issue is resolved and the first and second homology groups of any real flag manifold of type A are provided. In \cite{LR22c}, a novel characterization of the Bruhat order of the symmetric group is introduced via the permutation's Lehmer code. Building on the work developed in the aforementioned studies, we have advanced towards a comprehensive solution. Matszangosz \cite{Mat19}, through an alternative approach, also obtained a formula with certain combinatorial complexities. However, this method has the added advantage of being illustrated through the code diagram for permutations (see Figure \ref{fig:ex_sign}). Additionally, the third and fourth homology groups, along with the free part up to the sixth homology group, and their corresponding generators, are computed in full generality (see Proposition \ref{prop:betti_generators}, Theorems \ref{thm:H3},\ref{thm:H4}).

The paper is organized as follows. \Cref{sec:pre} introduces the primary definitions related to the combinatorics of the symmetric group, flag manifolds of type A, and their cellular decomposition. \Cref{sec:homology} presents the topological results concerning these manifolds. \Cref{sec:sign} details the process of deriving a formula for the sign of the boundary map coefficients. \Cref{sec:generators} focuses on the generators in the integral homology of flag manifolds of type A in specific low-dimensional cases.

\section{Preliminaries}\label{sec:pre}

Let $\N=\{1,2,3, \dots\}$ and $\Z$ be the set of integers. For $n,m\in\Z$, with $n\leqslant m$, denote the set $[n,m]=\{n,n+1, \dots, m\}$. For $n\in \N$, denote $[n]=[1,n]$.

\subsection{Combinatorics of \texorpdfstring{$S_n$}{Sn}}\label{subsec:prelim}

Let $\mathfrak{g}=\mathfrak{sl}_n(\mathbb{R})$ be the Lie algebra of type $A$ with real Lie group $G=\mathrm{Sl}(n,\mathbb{R})$. The Weyl group $\weyl$ of $\mathfrak{g}$ is the symmetric group $S_n$. The set of simple roots $\Sigma=\{a_{1},\dots, a_{n-1}\}$ is ordered as follows:

\begin{center}
\begin{tikzpicture} \dynkin[mark=o]{A}{};
\dynkinLabelRoot*{1}{a_1} \dynkinLabelRoot*{2}{a_2} \dynkinLabelRoot*{4}{a_{n-1}}
\end{tikzpicture}
\end{center}

As a Coxeter group, the permutation model for the symmetric group presents itself as a group generated by $s_i=(i,i+1)$, $i\in [n-1]$, with relations: (i) $s_i^2=1$, for $i\in [n-1]$; (ii) $s_is_j=s_js_i$, if $|i-j|>1$ for $i,j\in [n-1]$ and (iii) $s_is_{i+1}s_i=s_{i+1}s_is_{i+1}$, for $i \in [n-2]$. In this sense, we call $s_i=s_{a_i}$ a simple reflection for each $i\in [n-1]$, (ii) a commutation move, and (iii) a braid move. 

Each permutation $w \in S_n$ will be written in one-line notation $w=w(1)\cdots w(n)$.
We may consider the right and left actions of a simple reflection $s_i=(i,i+1)$: the right action swaps $w(i)$ and $w(i+1)$ (the values of $w$ at the \textit{positions} $i$ and $i+1$) while the left action exchanges the \textit{values} $i$ and $i+1$. A pair $(i,j)$ such that $i<j$ and $w(i)>w(j)$ is called an inversion of $w$.

We write $u\leq w$ if given a reduced decomposition $w=s_{j_{1}} \cdots s_{j_{r}}$ then $u=s_{j_{i_{1}}}\cdots s_{j_{i_{k}}}$, $1\leqslant i_1\leqslant \cdots \leqslant i_k\leqslant r$. It provides a partial order in $S_n$ called the strong Bruhat order.

Define the length $\ell(w)$ of $w \in \weyl$ as the number of simple reflections in any reduced decomposition of $w$ which is also the same as counting the number of its inversions. Also, define the code of $w$  by $\code(w)=\alpha=(\alpha_{1}, \dots, \alpha_{n-1})$ where $\alpha_{i} = \#\{k>i \ \colon w(k)<w(i)\}$ which gives the number of inversions to the right of each $w(i)$. Hence, $\ell(w)=\alpha_1+\cdots+\alpha_{n-1}$.

Define the set $\mathcal{C}_{n}$ as the Cartesian product $[0,n-1]\times [0,n-2] \times \cdots \times [0,1]$. Since $0\leqslant \alpha_{i}\leqslant n-i$, we have the code $\alpha \in \mathcal{C}_n$. If necessary, denote $\alpha_{n}=0$. The next proposition presents a useful relation between a descent of $w$ and its code.

\begin{proposition}[\cite{Mac91}, (1.24)]\label{prop:macdonald}
Let $i\in [n-1]$. Then $\alpha_{i} > \alpha_{i+1}$ if and only if $w(i) > w(i+1)$.
\end{proposition}

\subsection{\texorpdfstring{Quotients of $S_n$}{Quotients of Sn}}

For any $\Theta \subset \Sigma$, $\mathcal{W}_{\Theta} = \langle s_{i} \colon a_i \in \Theta\rangle $ is the standard parabolic subgroup. Denote by $\mathcal{W}^{\Theta} = \{ w \in \mathcal{W} \colon \ell(w) <\ell(ws_{i}), \forall a_i \in \Theta\}$ the set of representatives of $\mathcal{W}_{\Theta}$ of minimal length. 

For convenience, let us introduce a notation suitable for explicitly defining the root choice of $\Theta$. Let $\mathbf{k}=\mathbf{k}(\Theta)=\{k_{1},k_{2},\dots, k_{r}\}$ be the set of integers $0< k_1<k_2<\cdots<k_r < n$ such that $\Theta=\Sigma-\{a_{k_{1}},a_{k_{2}},\dots, a_{k_{r}}\}$, i.e., $\mathbf{k}(\Theta)$ are the indexes of all roots of $\Sigma$ that do not belong to $\Theta$.

In terms of the right action, the coset $w\weyl_{\Theta}$ consists of all $k_1!(k_2-k_1)!\cdots (n-k_r)!$ permutations obtained from $w$ by permuting the first $k_1$ entries among themselves, the following $k_2-k_1$ entries among themselves, and so on. Since the length $\ell(w)$ is exactly the number of inversions of $w$, the minimal representative inside a coset is given by the one for which these blocks consist of ordered entries. 

\begin{lemma}[\cite{Tan94}, Theorem A]\label{lem:minimal_representatives}
\begin{align}\label{eq:wthetamin}
\mathcal{W}^{\Theta} = \{ w=w(1)\cdots w(n) \in S_n \colon & w(1) < \cdots < w(k_1),  \nonumber \\
  & w(k_1+1) < \cdots <  w(k_2), \cdots, \\ & w(k_r+1)<\cdots<w(n)\} \nonumber.
\end{align}
\end{lemma}

That is, the elements of $\mathcal{W}^{\Theta}$ are precisely the permutations in $S_n$ for which the descents occur uniquely at the positions $k_1,k_2,\ldots, k_r$. It is immediate that 
\begin{equation*}
|\mathcal{W}^{\Theta}| = \dfrac{n!}{k_1!(k_2-k_1)!\cdots (n-k_r)!}.
\end{equation*}

In particular, when $\mathbf{k}=\{k\}$, the code applied to $\mathcal{W}^{\Theta}$ gives a bijection to the set $\mathcal{P}(k,n-k)$ of partitions whose Ferrers diagram fits in the rectangle $k\times(n-k)$, and Equation \eqref{eq:wthetamin} is written as
\begin{equation*}
\mathcal{W}^{\Theta} = \{ w=w(1)\cdots w(n)  \in S_n \colon w(1) < \cdots < w(k), w(k+1) < \cdots <\cdots<w(n)\}.
\end{equation*}

The idea of representing such permutations in terms of partitions may be generalized if instead the direct product of the set of partitions $\mathcal{P}(k_{1},n-k_1)\times \mathcal{P}(k_{2}-k_{1},n-k_{2})\times\cdots \times \mathcal{P}(k_{r}-k_{r-1},n-k_{r})$ is considered.

\begin{proposition}\label{prop:bijection_wk}
Given $\Theta$ such that $\mathbf{k}=\{k_{1},k_{2},\dots, k_{r}\}$ where $0< k_1<k_2<\cdots<k_r < n$, the code provides a bijection
\begin{equation*}
\weyl^{\Theta} \longleftrightarrow \mathcal{P}(k_{1},n-k_1)\times \mathcal{P}(k_{2}-k_{1},n-k_{2})\times\cdots \times \mathcal{P}(k_{r}-k_{r-1},n-k_{r}).
\end{equation*}
\end{proposition}
\begin{proof}
It is known that the code applied to any element of $S_{n}$ gives a bijection to the product $[0,n-1]\times [0,n-2] \times \cdots \times [0,1]$ (see \cite{LR22c}, Lemma 2.2). The proof follows directly from Proposition \ref{prop:macdonald} and Lemma \ref{lem:minimal_representatives}.
\end{proof}

Let $\alpha=(\alpha_{1},\dots,\alpha_{n-1})$ be the code of $w$. The proposition provides the decomposition of $\alpha$ into partitions $\alpha^{1},\dots,\alpha^{r}$ such that $\alpha^{j}\in \mathcal{P}(k_{j}-k_{j-1},n-k_{j})$ for each $j\in [r]$, i.e., it should be contained in a rectangle $(k_{j}-k_{j-1})\times (n-k_{j})$. 

For each permutation $w \in \weyl^{\Theta}$, define its corresponding diagram, which is called the code diagram of $w$. The code diagram of $w$ is the collection of left-justified boxes where the $i$-th row (counted from bottom to top) contains $\alpha_{i}$ boxes. Thus, $\alpha$  should be contained in the stack of rectangles of size $k_{1}\times (n-k_{1}), \dots, (k_{j}-k_{j-1})\times (n-k_{j}), \dots, (k_{r}-k_{r-1})\times (n-k_{r})$. Intuitively, the code diagram is a stack of minor Young diagrams in each rectangle.

\begin{figure}[hbtp]
\centering
\includegraphics[scale=1]{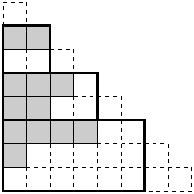}
\caption{The code diagram of $w=1\,3\, 7 \,5\,8\,2\,9\,4\,6$}
\label{fig:ex_diag}
\end{figure}

For example, let $n=9$ and $\Theta$ be such that $\mathbf{k}=(3,5,7)$. The permutation $w=1\,3\, 7 \,5\,8\,2\,9\,4\,6\in \weyl^{\Theta}$ corresponds to the code $\alpha=(0,1,4,2,3,0,2,0)$ whose code diagram is shown in Figure \ref{fig:ex_diag}.

\begin{remark}
Given $\Theta$ such that $\mathbf{k}=\{k_{1},k_{2},\dots, k_{r}\}$, there is an easy algorithm to retrieve $w\in\weyl^{\Theta}$ using the code diagram associated with the code $\alpha=(\alpha_1,\ldots, \alpha_{n-1})$. Given the Young diagram of $(\alpha_{1},\alpha_{2},\dots,\alpha_{k_{1}})$ inside the rectangle of size $k_{1}\times (n-k_{1})$, write a path along the boundary of the diagram starting from the SW corner to the NE corner of the rectangle containing the diagram. List the numbers $1,\dots, n$ for each step. Then, denote by $x_{1}(1)< x_{1}(2)< \cdots< x_{1}(k_{1})$ the numbers assigned to the vertical arrows and by $y_{1}(1)< y_{1}(2)< \cdots< y_{1}(n-k_{1})$ those of the horizontal arrows. 

Next, for the Young diagram of $(\alpha_{{k_1}+1},\dots,\alpha_{k_{2}})$ inside the rectangle of size $(k_{2}-k_{1})\times (n-k_{2})$, write a path along the boundary of the diagram and list the numbers $y_{1}$ for each arrow. As before, we denote by $x_{2}(1)< x_{2}(2)< \cdots< x_{2}(k_{2}-k_{1})$ the numbers assigned to the vertical arrows and by $y_{2}(1)< y_{2}(2)< \cdots< y_{2}(k_{2}-k_{1})$ those of the horizontal arrows.

Notice that $x_2$ is distinct from $x_1$. Repeat this process up to $(\alpha_{k_{r}+1},\dots,\alpha_{n-1})$. The corresponding permutation is the word obtained by joining $x_{1}$, $x_{2}$, $\ldots$, and $x_{r+1}$, where $x_{r+1}$ consists of the remaining numbers in the usual order.

For example, consider $n=9$, $\mathbf{k}=(3,5,7)$, and code $\alpha=(0,1,4,2,3,0,2,0)$. By the process described in Figure \ref{fig:constw}, $x_1=1<3<7$, $x_2=5<8$, $x_3=2<9$. Hence, $w=1\,3\, 7\, 5\,8\, 2\,9\, 4\,6$.

\begin{figure}[hbtp]
\centering
\includegraphics[scale=0.8]{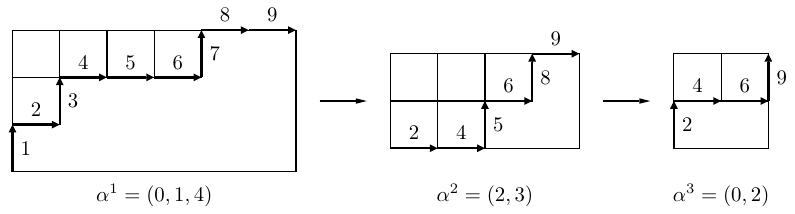}
\caption{Construction of $w$ from the code}
\label{fig:constw}
\end{figure}

\end{remark}

Arranging the boxes in a code diagram is advantageous because it provides an easy way to get the permutation $w$ in terms of simple reflections $s_{i}$. This is known as the \emph{row-reading expression} of $w \in \mathcal{W}^{\Theta}$. We begin to assign a simple reflection consecutively to each box from left to right and upwards, starting from $s_{1}$ in the bottom leftmost box in the staircase shape. Then, a reduced decomposition follows by reading each row in the diagram from right to left, and the rows from bottom to top. More specifically, define the reading for the $i$-th row by $\rr{w}_{i}= s_{\alpha_{i}+i-1} \cdot s_{\alpha_{i}+i-2} \cdots s_{i+1} \cdot s_{i}$ if $\alpha_{i}$ is non-zero, and $\rr{w}_{i} = e$ otherwise. The row-reading expression $\rr{w}$ of $w$ is
\begin{equation*}
\rr{w} = \rr{w}_{1}\cdots \rr{w}_{n-1}.
\end{equation*}

\begin{figure}[hbtp]
\centering
\includegraphics[scale=1]{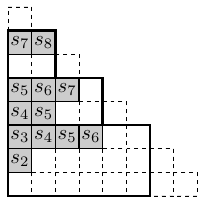}
\caption{Row-reading of $w=1\,3\, 7 \,5\,8\,2\,9\,4\,6$}
\label{fig:rowreading}
\end{figure}

The row-reading for the permutation $w=1\,3\, 7 \,5\,8\,2\,9\,4\,6$ is illustrated in Figure \ref{fig:rowreading}. Then, $\rr{w}=s_{2}\cdot s_{6}s_{5}s_{4}s_{3}  \cdot s_{5}s_{4} \cdot s_{7}s_{6}s_{5} \cdot s_{8}s_{7}$.

\begin{remark}

The code diagram of any permutation $w\in \weyl=S_n$ should fit in a staircase shape since it corresponds to the choice of $\Theta=\emptyset$, i.e., $\mathbf{k}=[n-1]$ (see \cite{LR22c} for details).
\end{remark}

\subsection{Flags of \texorpdfstring{$\mathrm{Sl}_n(\mathbb{R})$}{Sln(R)}}\label{subsec:flags}

Flag manifolds are defined as homogeneous spaces of type $G/P$ where $G$ is a non-compact semi-simple Lie group and $P$ is a parabolic subgroup of $G$. The present work focuses on the flag manifolds of type A, i.e., the flag manifolds of $G=\mathrm{Sl}_n(\mathbb{R})$ with the corresponding real semi-simple Lie algebra $\mathfrak{sl}_n(\mathbb{R})$ which is a split real form of $\mathfrak{sl}_n(\mathbb{C})$. The Iwasawa decomposition provides $G=KAN$ where $K=\mathrm{SO}(n)$ is the (compact) group of orthogonal matrices, $A$ of diagonal matrices, and $N$ given by the unipotent matrices. The corresponding decomposition for its Lie algebra is $\mathfrak{sl}_n(\mathbb{R})=\mathfrak{so}(n)\oplus \mathfrak{a}\oplus \mathfrak{n}$, where $\mathfrak{a}$ is the subalgebra of diagonal matrices with zero trace and $\mathfrak{n}$ is the subalgebra of the strict upper triangular matrices (for details, see \cite{HN12}, Chapter 13). In this context, the simple roots of $\Sigma=\{a_1,\ldots, a_{n-1}\}$ are given by $a_i = \varepsilon_{i}-\varepsilon_{i+1}$, where $\varepsilon_i \in \mathfrak{a}^{\ast}$ is defined by $\varepsilon_i(\mathrm{diag}(\lambda_1,\ldots, \lambda_n))=\lambda_i$.

For any $\Theta\subset \Sigma$, the corresponding flag manifold $\mathbb{F}_{\Theta}=\mathrm{Sl}_n(\mathbb{R})/P_{\Theta}$, where $P_{\Theta}$ is a standard parabolic subgroup containing $P$. If $\Theta=\emptyset$, omit the subscript $\Theta$ and call $\mathbb{F}$ a maximal flag. Otherwise, $\mathbb{F}_{\Theta}$ is called a partial flag manifold.

In the context of type $A$, each $\mathbb{F}_\Theta$ is the manifold of a nested sequence of subspaces of $\mathbb{R}^n$, while the parabolic subgroup has a block structure related to the roots of $\Theta$. Explicitly, accordingly to the notation introduced in \Cref{subsec:prelim},  if $\mathbf{k}=\{k_1,\dots,k_r\}$ then
\begin{equation*}
\mathbb{F}_{\Theta}=\mathbb{F}(k_1,k_2-k_1,\dots,n-k_r)=\{(V_{1}\subset\cdots \subset V_{r}\subset \mathbb{R}^n)\colon \dim V_{i}=k_i-k_{i-1}\}.
\end{equation*}

The Bruhat decomposition of any $\mathbb{F}_{\Theta}=G/P_{\Theta}$ presents the flag as the union of $N$-orbits parametrized by the minimal representatives $\weyl^{\Theta}$ of $\weyl_{\Theta}$ in $\weyl$, i.e., 
\begin{equation*}
\mathbb{F}_{\Theta} = \bigsqcup_{w\in \weyl^{\Theta}} N\cdot wb_{\Theta},
\end{equation*}
where $b_{\Theta}=1\cdot P_{\Theta}$ is the base point. For each $w \in \weyl^{\Theta}$, $N\cdot wb_{\Theta}$ is called a Bruhat cell. The closure of a Bruhat cell is called a Schubert cell and is denoted by $\mathcal{S}_w^{\Theta}$.

The Schubert cells provide a cellular structure to the flag manifolds with 
\begin{equation*}
    \mathcal{S}^{\Theta}_w=\bigcup_{u\leq w} S^{\Theta}_u, 
\end{equation*}
where $u \in \weyl^{\Theta}$. Since $\mathfrak{g}$ is a split real form, then $\dim \mathcal{S}_w^{\Theta}=\ell(w)$, for each $w \in \weyl^{\Theta}$.

The natural projection map $\pi_{\Theta}: \mathbb{F}\rightarrow \mathbb{F}_{\Theta}$ is equivariant with respect to the $G$-action and, more generally, if $\Delta \subset \Theta \subset \Sigma$ there is also the projection map $\pi^{\Delta}_{\Theta}: \mathbb{F}_{\Delta}\rightarrow \mathbb{F}_{\Theta}$. Notice that several Schubert cells $S_w$ in the maximal flag manifold project into a Schubert cell $\mathcal{S}_w^{\Theta}$ in the partial flag manifold and there is a unique Schubert cell among them with the same dimension (see \cite{RSm19}, Lemma 3.1).  

It is possible to obtain explicit parametrizations for Schubert cells inside $K$, the compact group given by the Iwasawa decomposition, since there is a diffeomorphism $\mathbb{F}_{\Theta}\cong K/K_{\Theta}$, where $K_{\Theta}=K\cap P_{\Theta}$.

Let us explore such a parametrization for flags of type $A$. It will be sufficient to present such parametrizations in the maximal flag manifold of $\mathrm{Sl}_n(\mathbb{R})$. Denote by $[0,\pi]$ the interval of real numbers $0\leq t\leq \pi$, and by $E_{k,l}$ the $n\times n$ integer matrix such that it is $1$ at position $(k,l)$ and $0$ otherwise. For each simple reflection $s_{i}=s_{a_{i}} \in \Sigma$, consider the matrix $A_{s_i}=A_{i}$ defined by $A_{i} = E_{i,i+1} - E_{i+1,i}$, for $i\in[n-1]$. By Proposition 1.3 of \cite{RSm19}, if $w=s_1\cdots s_{\ell}$ is a reduced decomposition, then the corresponding Schubert cell is parametrized by $(t_1, \ldots, t_{\ell}) \in [0,\pi]^{\ell} \mapsto e^{t_1A_1} \cdots e^{t_{\ell}A_{\ell}}\cdot b_0$, where the exponential of the matrix $t_iA_{i}$, for $t_i\in[0,\pi]$, is the rotation matrix.

For the type A flag manifold ${\mathbb{F}}_{\Theta}$, by the permutation model for $S_n$, a Schubert cell is uniquely defined by a permutation whose dimension is the number of boxes in the corresponding diagram according to Proposition \ref{prop:bijection_wk}. Let $\pi_{\Theta}:\mathbb{F}\rightarrow \mathbb{F}_{\Theta}$ be the natural projection. The image of $\mathcal{S}_w$ through $\pi_{\Theta}$, for a permutation $w\in S_n$, is the Schubert cell $\mathcal{S}^{\Theta}_{v}$ for which $v\in \weyl^{\Theta}$, i.e., $v$ is the permutation obtained by ordering the entries of $w$ according to the partition $0<k_1<\cdots<k_r$ such that the descents occur at positions $k_1, \ldots, k_r$. 

\section{Homology of flag manifolds of type A}\label{sec:homology}
This section is devoted to establish the formula for the coefficients of the boundary map of real flag manifolds of type A as obtained in \cite{LR22a}, which is not definitive for the purposes of this paper. It remains to solve the problem of the ambiguity of the signs, which will be addressed in the next Section.

Some results on the Poincaré polynomial that has recently been described by He \cite{He2017, He2020} are also stated with a compatible notation. 

\subsection{Coefficients of the boundary map}

Let us start briefly in the context of the maximal flag manifold $\mathbb{F}$. For every $w\in S_{n}$, fix the reduced decomposition of $w$ the row-reading $\rr{w}=s_{1} \cdots s_{\ell(w)}$, with $\ell=\ell(w)$.
Let $\mathcal{C}$ be the $\mathbb{Z}$-module freely generated by all Schubert cells $\mathcal{S}_{w}$, $w\in S_n$. The boundary map $\partial$ defined over $\mathcal{C}$ is given by $\partial \mathcal{S}_{w}=\sum_{w^{\prime}}c(w,w^{\prime })\mathcal{S}_{w^{\prime }}$, where $c(w,w^{\prime })\in \mathbb{Z}$. The non-trivial coefficients must be equal to $\pm 2$ (\cite{RSm19}, Theorem 2.2) and they will occur only in the cases where $w$ covers $w'$, i.e., if $w=s_{1}\cdots s_{\ell}$ is a reduced decomposition of $w\in \mathcal{W}$, then there exists a unique index $I$ (which depends on the pair $w$ and $w'$) such that $w'= s_{1}\cdots \widehat{s_{\I}}\cdots s_{\ell}$ is a uniquely defined reduced decomposition (\cite{RSm19}, Proposition 1.10). Denote this reduced decomposition of $w'$ by $\rrt{w}{I}$, which may be different from its row-reading expression $\rr{w}'$.

The following result characterizes the covering relations for the Bruhat order of $S_n$.

\begin{lemma}[\cite{BB05}, Lemma 2.1.4]\label{lem:bjorner}
    Let $w,w'\in S_{n}$. Then, $w$ covers $w'$ in the Bruhat order if and only if $w = w' \cdot (i,j)$ for some transposition $(i,j)$ with $i<j$ such that $w'(i) < w'(j)$ and there does not exist any $k$ such that $i<k<j$, $w'(i) < w'(k) <w'(j)$.
\end{lemma}

The lemma says that if $w=w(1)\cdots w(n)$ is the one-line notation of $w\in S_{n}$ then $w'$ is covered by $w$ if and only if the one-line notation of $w'$ is obtained from $w$ by switching the values in positions $i$ and $j$, for some pair $i<j$ and such that no value between positions $i$ and $j$ lies in $[w(j),w(i)]$.

According to \cite{RSm19} Proposition 1.9, there are characteristic maps for $\mathcal{S}_{w'}$ given by $\Phi_{\rr{w}'}\colon B_{\rr{w}'}\rightarrow \mathcal{S}_{w'}$ and $\Phi_{\rrt{w}{I}}\colon B_{\rrt{w}{I}}\rightarrow \mathcal{S}_{w'}$, where $B_{\rr{w}'}$ and $B_{\rrt{w}{I}}$ are balls of dimension $\ell(w')$. Different reduced decompositions of $w'$ may provide distinct characteristic maps for $\mathcal{S}_{w'}$ since $\Phi_{{\rr{w}}'}$ and $\Phi_{\rrt{w}{I}}$ are, respectively, obtained from the row-reading of $w'$ and from the fixed choice of a reduced decomposition for $w$.

\begin{theorem}\label{thm:rabelosanmartin}
Let $w,w'\in S_{n}$ be a pair such that $w$ covers $w'$. Then,
\begin{equation}\label{eq:coefficient}
    c(w,w^{\prime})=(-1)^{\I}\cdot\deg(\Phi_{\rrt{w}{I}}^{-1} \circ \Phi_{\rr{w'}})\cdot(1+(-1)^{j-i}),
\end{equation}
where $I$ is the index such that row-reading $\rr{w}=s_{1}\cdots s_{\I}\cdots s_{\ell}$ provides the reduced decomposition $\rrt{w}{I}=s_{1}\cdots \widehat{s_{\I}}\cdots s_{\ell}$ of $w'$, and the pair $i,j$ satisfies $w = w' \cdot (i,j)$ as in \Cref{lem:bjorner}.
\end{theorem}
\begin{proof}
This follows directly from Theorem 2.8 in \cite{RSm19} and Proposition 5.2 of \cite{LR22a}.
\end{proof}

\begin{remark}\label{rem:degree}
Denote by $\sigma_{w'}=\mathcal{S}_{w'}/(\mathcal{S}_{w'}\setminus N\cdot w'b_0)$ the space given by identifying the complement of the Bruhat cell $N\cdot w'b_0$ to a point. The characteristic maps $\Phi_{\rr{w}'}$ and $\Phi_{\rrt{w}{I}}$ induce homeomorphims $S^{\ell(w')}\to \sigma_{w'}$ by collapsing the boundary of the corresponding balls to points. Hence $\Phi_{w'}^{-1}\circ \Phi_{\widehat{w}_\I}$ is a map between spheres where $\Phi_{w'}^{-1}$ is the inverse of the corresponding induced homeomorphism $S^{\ell(w')}\to \sigma_{w'}$.
\begin{center}
$
\xymatrix{
B_{\rr{w}'} \ar[r]^{\Phi_{\rr{w}'}} \ar[d] & \mathcal{S}_{w'} \ar[d] & B_{\rrt{w}{I}} \ar[l]_{\Phi_{\rrt{w}{I}}}\ar[d] \\
S^{\ell(w')}=B_{\rr{w}'}/\partial(B_{\rr{w}'}) \ar[r]  & \sigma_{w'} & B_{\rrt{w}{I}}/ \partial(B_{\rrt{w}{I}})=S^{\ell(w')} \ar[l]
}
$
\end{center}
\end{remark}

\begin{remark} If both reduced decompositions $\rr{w}'$ and $\rrt{w}{I}$ coincide then $\deg(\Phi_{\rrt{w}{I}}^{-1} \circ \Phi_{\rr{w'}})=1$.
\end{remark}

Consider the partial flag manifold $\mathbb{F}_{\Theta}$.  Let $\mathcal{C}^{\Theta}$ be the $\mathbb{Z}$-module freely generated by $\mathcal{S}_{w}^{\Theta}$, for every element $w$ of $\weyl^{\Theta}$. 

Recall that $\pi_{\Theta}:\mathbb{F}\to\mathbb{F}_{\Theta}$ projects Schubert cells $\mathcal{S}_{w}$ into $\mathcal{S}_{w}^{\Theta}$. The integral homology of the flag manifold $\mathbb{F}_{\Theta}$ is isomorphic to the homology of $(\mathcal{C}^{\Theta},\partial^{\Theta})$, where $\partial^{\Theta}$ is obtained by restricting $\partial$ and projecting it onto $\mathcal{C}^{\Theta}$. If we write the boundary map $\partial^{\Theta} :\mathcal{C}^{\Theta} \rightarrow \mathcal{C}^{\Theta}$ as 
\begin{equation*}
\partial^{\Theta} \mathcal{S}_{w}^{\Theta}=\sum_{w'\in \weyl^{\Theta}}c^{\Theta}(w,w')\mathcal{S}_{w'}^{\Theta}
\end{equation*}
then, by \cite{RSm19}, Theorem 3.4, it follows that $c^{\Theta}(w,w')=c(w,w')$ for $w,w'\in\weyl^{\Theta}$, where $c(w,w')$ is the coefficient of $\partial$ in the maximal flag manifold.

\subsection{Poincaré polynomial}

For a flag manifold $\mathbb{F}_{\Theta}$, write the integer and $\Z_2$ homology groups as
\begin{align*}
    H_i(\mathbb{F}_{\Theta}, \Z) &\cong \Z^{\beta_i} \oplus (\Z_2)^{T_i}\\
    H_i(\mathbb{F}_{\Theta}, \Z_2) &\cong (\Z_2)^{B_i}
\end{align*}
where $\beta_i$ are the Betti numbers, $(\Z_2)^{T_i}$ is the torsion, and $B_i$ counts the number of Schubert cells of dimension $i$ (see \cite{RSm19}, Corollary 2.3). 

Following \cite{He2020}, denote the Poincaré polynomial in integer and mod-2 coefficients, and the generating polynomial of the  torsions in integral homology, respectively, by
\begin{align*}
    FP(t)&=\sum_{i=0}^{n} \beta_i t^i; & P(t) & = \sum_{i=0}^{n} B_i t^i; & TP(t) & = \sum_{i=0}^{n} T_i t^i.
\end{align*}

Given $n\in\N$, the $t$-analogue of $n$ and $n!$ are, respectively, the polynomials $(n)_t = 1 + t + \cdots + t^{n-1}=\frac{1-t^n}{1-t}$ and $(n)_t!=(1)_t(2)_t\cdots (n)_t = (1+t)(1+t+t^2)\cdots (1+t+ \cdots + t^{n-1})$.

Given $n_1,\dots,n_k$ positive integers such that $n=n_1+\cdots+n_k$, define the $t$-multinomial coefficient by
\begin{equation*}
    \binom{n}{n_1,\dots,n_k}_t = \frac{(n)_t!}{(n_1)_t!\cdots (n_k)_t!}.
\end{equation*}

Let $\Theta\subset \Sigma$. Denote by $\Theta'= \{a_{k_{1}},\dots, a_{k_{s}}\}$ the complement of $\Theta$, with $\mathbf{k}=\{k_1,\ldots, k_s\}$ as defined in Section \ref{subsec:prelim}. Define $L = L(\Theta)$ by
\begin{equation*}
L = \left\lfloor \dfrac{k_{1}}{2} \right\rfloor + \left\lfloor \dfrac{k_{2}-k_{1}}{2} \right\rfloor + \cdots + \left\lfloor \dfrac{k_{s}-k_{s-1}}{2} \right\rfloor + \left\lfloor \dfrac{n-k_{s}}{2} \right\rfloor,
\end{equation*}
where $\lfloor \cdot \rfloor$ is the floor function. Clearly, $L(\emptyset) = 0$.

The next proposition describes the Poincaré polynomials for any real flag manifolds of type A. 

\begin{proposition}[\cite{Bo53}, and \cite{He2020} Prop. 3.14]\label{prop:poincare}
Let $\Theta\subset \Sigma$.
Then, the mod-2 Poincaré polynomial of the partial real flag manifold $\mathbb{F}_{\Theta}$ of type A is
\begin{equation*}
    P(t) = \binom{n}{k_1, k_2-k_1, \dots, n-k_s}_{t},
\end{equation*}
and the Poincaré polynomial in integer coefficients is
\begin{itemize}
\item If $n$ is odd or $L=\frac{n}{2}$:
\begin{equation*}
FP(t) = \binom{L}{ \left\lfloor \frac{k_{1}}{2} \right\rfloor , \left\lfloor \frac{k_{2}-k_{1}}{2} \right\rfloor , \cdots , \left\lfloor \frac{n-k_{s}}{2} \right\rfloor}_{t^{4}}\cdot \prod_{i = L}^{\left\lfloor \frac{n-1}{2}\right\rfloor-1} (1+t^{4i+3})
\end{equation*}
\item If $n$ is even and $L\neq \frac{n}{2}$:
\begin{equation*}
FP(t) = (1+t^{n-1})\binom{L}{ \left\lfloor \frac{k_{1}}{2} \right\rfloor , \left\lfloor \frac{k_{2}-k_{1}}{2} \right\rfloor , \cdots , \left\lfloor \frac{n-k_{s}}{2} \right\rfloor}_{t^{4}}\cdot \prod_{i = L}^{\left\lfloor \frac{n-1}{2}\right\rfloor-1} (1+t^{4i+3})
\end{equation*}
\end{itemize}
\end{proposition}

\begin{remark} In \cite{He2020}, the cohomology ring is determined using $\mathbb{Z}[1/2]$ coefficients. However, the Universal Coefficient Theorem ensures that the Poincaré polynomial for $\mathbb{Z}[1/2]$ coefficients coincides with the polynomial for integer cohomology. In essence, $\mathbb{Z}[1/2]$ cohomology removes the torsion while preserving the free part. 
\end{remark}

\begin{lemma}[\cite{He2017}, Lemma 3.2] \label{lem:chenhe} The generating polynomial of the torsion in the integral homology for the partial flag real flag manifolds $\mathbb{F}_{\Theta}$ of type A is
\begin{equation*}
    TP(t) = \frac{P(t)-FP(t)}{1+t}.
\end{equation*}    
\end{lemma} 

For groups of type A, it is known that if $\Theta_{1},\Theta_{2} \subset \Sigma$ and the corresponding Dynkin diagrams are the same, then the flag manifolds are homeomorphic. In fact, for $\Theta'_1 = \{a_{k_{1}},\dots, a_{k_{s}}\}$ and $\Theta'_2 = \{a_{l_{1}},\dots, a_{l_{s}}\}$ each can be written as $\mathbb{F}_{\Theta_1} \cong O(n)/O(k_1)\times O(k_2-k_1) \times \cdots \times O(n-k_s)$ and $\mathbb{F}_{\Theta_2} \cong O(n)/O(l_1)\times O(l_2-l_1) \times \cdots \times O(n-l_s)$. Hence, both are homeomorphic since $\{k_{1},k_2-k_1,\dots, n-k_{s}\} = \{l_{1},l_2-l_1,\dots, n-l_{s}\}$ as multisets. In this case, $H_i(\mathbb{F}_{\Theta_{1}}, \Z)\cong H_i(\mathbb{F}_{\Theta_{2}}, \Z)$ for any $i$. Therefore, the Betti numbers $\beta_i$ and the torsion coefficients $T_i$ depend only on the cardinalities of the connected components of the Dynkin diagram of $\Theta$.

\section{Sign of the boundary map}\label{sec:sign}

In the context of flags of type A, it is possible to go even further than simply saying that the sign of the boundary map is $0$ or $\pm 2$ by constructing a complete algorithm to compute the coefficients of the boundary map according to Theorem \ref{thm:rabelosanmartin}. Indeed, it remains to obtain the degree of $\Phi_{\rrt{w}{I}}^{-1} \circ \Phi_{\rr{w'}}$ in Equation \eqref{eq:coefficient}. This next step is strongly based on the combinatorics of the codes introduced in \cite{LR22c}. As mentioned in the Introduction, this process is not unique in the sense that it is possible to make another choice of signs for $c(w,w')$ such as, for example, that given by \cite{Mat19}.

\subsection{Another criteria for covering relations}

Let $w\in S_{n}$ and $\alpha=\code(w)$. Given $i\in[n]$ and $j \in \N$, define the coefficients $M_{i,j}$ for $w$ recursively as follows:
\begin{itemize}
\item If $j\leqslant i+1$, then $M_{i,j}(w) = 0$;
\item If $j>i+1$, then $M_{i,j}(w) = M_{i,j-1}(w) + \left\{
\begin{array}{cl}
1, & \mbox{if } \alpha_{j-1} < \alpha_{i} - M_{i,j-1}(w);\\
0, & \mbox{if } \alpha_{j-1} \geqslant \alpha_{i} - M_{i,j-1}(w).
\end{array}\right.$
\end{itemize}

It follows that $M_{i,k+1}(w)-M_{i,k}(w)$ is either zero or one. These coefficients are built using only the code of $w$. The next proposition shows that it is also possible to define them using the permutation $w$.

\begin{proposition}[\cite{LR22c}, Lemma 4.1]
    For $1\leqslant i <j \leqslant n$, 
    \begin{equation}
        M_{i,j}(w) = \#\{k \colon i<k<j \mbox{ and } w(k)<w(i) \}.\label{eq:cw}
    \end{equation}
\end{proposition}

Denoncourt in \cite{Den13} calls the $n\times n$ matrix $M(w)_{i,j}=(M_{i,j}(w))$ the extend matrix of $w$. Notice that $M_{i,i+1}(w)=0$ and $M_{i,n+1}(w) = \alpha_{i}$ for $i\in[n]$.

The coefficients $M_{i,j}$ are useful in describing the Bruhat order. From now on, given $w,w'\in S_{n}$, we will denote $\alpha = \code(w)$ and $\alpha' = \code(w')$.

\begin{theorem}[\cite{LR22c}, Theorem 4.4]\label{thm:codecovering2}
    Let $w,w'\in S_{n}$. Then, $w$ covers $w'$ with $w'=w\cdot (i,j)$, $1\leq i<j\leq n$, if, and only if, the following conditions are satisfied:
    \begin{align}
        \alpha_{i}' &\leqslant \alpha_{i} - 1;\label{eq:a1}\tag{a1}\\
        \alpha_{j}' & =  \alpha_{j} + \alpha_{i}  -\alpha_{i}'  - 1;\label{eq:a2}\tag{a2}\\
        \alpha_{k}' & = \alpha_{k} \mbox{ for every } k \neq i \mbox{ and } k \neq j;\label{eq:a3}\tag{a3}\\
        M_{i,j}(w) & = M_{i,j}(w') = \alpha_{i}'- \alpha_{j}.\label{eq:a4}\tag{a4}
    \end{align}
\end{theorem}

For instance, let $w,w' \in S_7$ be the pair of permutations given by $w=4\, 6\, 7\, 2\, 3\, 1\, 5$ and $w'=4\, 2\, 7\, 6\, 3\, 1\, 5$, i.e., $w=w'\cdot (2,4)$. Hence $w,w'$ is a covering pair by \Cref{lem:bjorner}. The corresponding codes $\alpha = (3,4,4,1,1,0)$ and $\alpha' = (3,1,4,3,1,0)$ clearly satisfy \eqref{eq:a1},\eqref{eq:a2} and \eqref{eq:a3}. Since $M_{2,4}(w)=M_{2,4}(w')=0=\alpha_2'-\alpha_4$, \eqref{eq:a4} is also true, giving another proof of the covering relation. For the purpose of this verification, the criteria of \Cref{lem:bjorner} is more straightforward than those given by \Cref{thm:codecovering2}. However, the former provides more precise information about the reduced decompositions involved. This will become clearer in the next section.

\subsection{Counting moves between reduced decompositions of \texorpdfstring{$w^{\prime}$}{w'}}

For a given element $w \in \mathcal{W}$, define the graph $R(w)$ on the set of reduced words for $w$. The vertices of $R(w)$ represent the reduced expressions of $w$ in terms of simple reflections. An edge connects two vertices if and only if their corresponding reduced words are related by either a single commutation or a braid move. According to Theorem 3.3.1(ii) in \cite{BB05}, any reduced decomposition of $w$ can be transformed into any other by a sequence of commutation and braid moves. Therefore, the graph $R(w)$ is connected.

Given $w, w' \in S_n$ such that $w$ covers $w'$, it is possible that $\rrt{w}{I}$ and $\rr{w'}$ are distinct reduced decompositions of $w'$. To address this issue, construct a sequence of commutations and braid moves that transforms $\rrt{w}{I}$ into $\rr{w'}$. This sequence corresponds to a specific path in the graph $R(w) $. 

\begin{proposition}[\cite{BCL15}, Corollary 5] \label{prop:loopingw} Let $\mathcal{W}$ be a finite Coxeter group and $w \in \mathcal{W}$. Then any loop in the graph $R(w)$ contains an even number of edges. 
\end{proposition}

\begin{proposition}[\cite{LR22c}, Proposition 4.3]\label{prop:moves_property2}
    Suppose that $w$ covers $w'$ with $w'=w\cdot (i,j)$, $1\leq i<j\leq n$, such that $\rrt{w}{I}$ differs from $\rr{w}'$.
    Then, for every $k\in [i+1,j-1]$ and $m\in [\alpha_{i}'+1,\alpha_{i}-1]$,
    \begin{equation}
        s_{m+i} \cdot \rr{w}_{i}' = \rr{w}_{i}' \cdot s_{m+i}\label{eq:lemaproperties1},
    \end{equation} 
    where the reduced decomposition on the left produces the one on the right after a sequence of $\alpha_i'$ commutations and no braid move, and
    \begin{equation}
        s_{m+k-1 - M_{i,k}(w)} \cdot \rr{w}_{k} = \rr{w}_{k} \cdot s_{m+k - M_{i,k+1}(w)}\label{eq:lemaproperties2},
    \end{equation}
    where the reduced decomposition on the left produces the one on the right after a sequence of $1-M_{i,k+1}(w)+M_{i,k}(w)$ braid moves and $\alpha_k - 2(1-M_{i,k+1}(w)+M_{i,k}(w))$ commutations.
\end{proposition}

\begin{figure}[hbtp]
    \centering
    \includegraphics[scale=0.7]{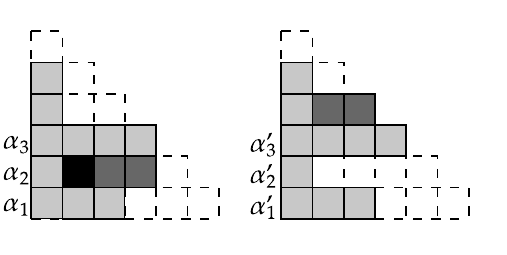}
    \caption{Diagram of $w=4\, 6\, 7\, 2\, 3\, 1\, 5$ in the left and $w'=4\, 2\, 7\, 6\, 3\, 1\, 5$ in the right.}
    \label{fig:ex_sign}
\end{figure}

For example, let $w,w' \in S_7$ be the covering pair of permutations given by $w=4\, 6\, 7\, 2\, 3\, 1\, 5$ and $w'=4\, 2\, 7\, 6\, 3\, 1\, 5$, as previously. \Cref{fig:ex_sign} shows the code diagrams of $w$ and $w'$. The black box is the one which is removed to obtain $w'$. Recall that $w=w'\cdot (2,4)$ such that the corresponding diagrams differ exactly in the second and fourth lines. 
\Cref{prop:moves_property2} explains how one can move the remaining two dark gray boxes of the second line at right of the removed box to the fourth line. By \Cref{eq:lemaproperties1} applied twice, it follows that
\begin{equation*}
s_5s_4\cdot \rr{w}'_2 = \rr{w}'_2\cdot s_5s_4. 
\end{equation*}
By \Cref{eq:lemaproperties2} applied twice, it follows that
\begin{equation*}
s_5s_4\cdot \rr{w}'_3 = \rr{w}'_3\cdot s_6s_5. 
\end{equation*}
These equations explain why one can obtain $\rr{w}'$ from $\rrt{w}{I}$ by plugging them into the reduced decomposition of $w$:
\begin{align*}
\rrt{w}{I} &= \rr{w}'_1\cdot (s_5s_4\widehat{s_3}s_2) \cdot \rr{w}'_3 \cdot \rr{w}_4 \cdot \rr{w}'_5 \\
&=\rr{w}'_1\cdot s_2 \cdot (s_5s_4 \rr{w}'_3) \cdot \rr{w}_2 \cdot \rr{w}'_5 \\
&=\rr{w}'_1\cdot s_2\cdot \rr{w}'_3 \cdot (s_6s_5 \rr{w}_2) \cdot \rr{w}'_5 = \rr{w}'.
\end{align*}

\begin{remark} This process is also called a lader move in the diagram. For details, see \cite{LR22c}, Section 3.2.
\end{remark}

\begin{proposition}\label{prop:countbraid}
    Let $w,w'\in S_{n}$ be such that $w$ covers $w'$ with $w'=w\cdot (i,j)$, $1\leq i<j\leq n$, and $\rrt{w}{I}$ differs from $\rr{w}'$. There exists a sequence $f_1 f_2 \cdots f_p$ of either commutations or braid moves such that the number of commutations and braid moves required to transform $\rrt{w}{I}$ into $\rr{w'}$ is
    \begin{align*}
    \#\{k\colon f_k \mbox{ is a braid move}\} &= (\alpha_{i}-\alpha_{i}' - 1) (j-i-1+\alpha_{j}-\alpha_{i}'),\\
    \#\{k\colon f_k \mbox{ is a commutation}\} &= (\alpha_{i}-\alpha_{i}' - 1) \left(\left(\sum_{k=i}^{j-1} \alpha_{k}'\right) - 2 (j-i-1+\alpha_{j}-\alpha_{i}')\right).
    \end{align*}

    Moreover, given any other sequence of moves $g_1 g_2 \cdots g_q$ from $\rrt{w}{I}$ to $\rr{w'}$, then the number of braid moves (resp. commutations) in $g_k$ and $f_k$ have the same parity.
\end{proposition}
\begin{proof}

Since $\rrt{w}{I}$ differs from $\rr{w}'$, it follows that $\alpha_{i} > \alpha_{i}'+1$. As a consequence of conditions \eqref{eq:a1}, \eqref{eq:a2}, and \eqref{eq:a3} of \Cref{thm:codecovering2} we have that
\begin{align}
\rr{w}_{i}&=s_{\alpha_{i}+i-1} s_{\alpha_{i}+i-2} \cdots s_{\alpha_{i}'+i} \cdot \rr{w}_{i}';\label{eq:rri} \\
\rr{w}_{j}'&= s_{\alpha_{j}'+j-1} s_{\alpha_{j}'+j-2} \cdots s_{\alpha_{j}+j} \cdot \rr{w}_{j};\label{eq:rrj}\\
\rr{w}_{k}&=\rr{w}_{k}', \mbox{ for $k\neq i$ and $k\neq j$.} \label{eq:rrk}
\end{align}

Moreover, condition \eqref{eq:a2} also says that
\begin{equation*}
\ell(s_{\alpha_{j}'+j-1} s_{\alpha_{j}'+j-2} \cdots s_{\alpha_{j}+j}) = \ell(s_{\alpha_{i}+i-1} s_{\alpha_{i}+i-2} \cdots s_{\alpha_{i}'+i+1}).
\end{equation*}

The reduced decomposition $\rrt{w}{I}$ corresponds to removing the simple reflection $s_{\alpha_i'+i}$ from $\rr{w}$. Notice that Equation \eqref{eq:rri} shows the relationship between the row reading of $\rr{w}_i$ and $\rr{w}'_i$. It follows that $\rrt{w}{I}$ is obtained from $\rr{w}$ by deleting the first box in the $i$-row just after the word $\rr{w}'_i$ which reads $s_{\alpha_i'+i}$. Then, $\rrt{w}{I}$ and $\rr{w}'$ are written as follows:
\begin{align*}
\rrt{w}{I} &= \rr{w}_{1} \cdots \rr{w}_{i-1} \cdot ( s_{\alpha_{i}+i-1} s_{\alpha_{i}+i-2} \cdots s_{\alpha_{i}'+i+1} )\cdot \rr{w}_{i}' \cdot \rr{w}_{i+1}  \cdots \rr{w}_{n-1},\\
\rr{w'} &=\rr{w}_{1}\cdots \rr{w}_{i-1} \cdot \rr{w}_{i}'  \cdots \rr{w}_{j-1}' \cdot ( s_{\alpha_{j}'+j-1} s_{\alpha_{j}'+j-2} \cdots s_{\alpha_{j}+j} ) \cdot \rr{w}_{j} \cdots \rr{w}_{n-1}.
\end{align*}

Now, let us describe how to transform $\rrt{w}{I}$ into $\rr{w'}$ by exhibiting a sequence of moves $f_1 f_2 \cdots f_p$ such that
\begin{equation*}
\rrt{w}{I} = \rr{w}^{0} \xrightarrow{f_1} \rr{w}^{1} \xrightarrow{f_2} \rr{w}^{2} \xrightarrow{f_3} \cdots \xrightarrow{f_p} \rr{w}^{r} = \rr{w'}.
\end{equation*}

For each $m\in [\alpha_{i}'+1,\alpha_{i}-1]$, we can apply \Cref{prop:moves_property2} iteratively to get the sequence of moves as follows:
\begin{align*}
s_{m+i} \cdot \rr{w}_{i}' \cdots \rr{w}_{j-1}' &= \rr{w}_{i}' \cdot s_{m+i} \cdot \rr{w}_{i+1}' \cdot \rr{w}_{i+2}' \cdots \rr{w}_{j-1}' \nonumber \\
&= \rr{w}_{i}'  \cdot \rr{w}_{i+1}' \cdot s_{m+i+1-M_{i,i+2}(w)} \cdot \rr{w}_{i+2}' \cdots \rr{w}_{j-1}' \nonumber \\
&= \rr{w}_{i}'  \cdot \rr{w}_{i+1}' \cdot \rr{w}_{i+2}' \cdot s_{m+i+2-M_{i,i+3}(w)}  \cdots \rr{w}_{j-1}' \label{eq:commutestep}\\
&\qquad  \vdots \nonumber\\
& = \rr{w}_{i}'  \cdots \rr{w}_{j-1}' \cdot s_{m+j-1- M_{i,j}(w)} = \rr{w}_{i}'  \cdots \rr{w}_{j-1}' \cdot s_{m+j-1+\alpha_{j}-\alpha_{i}'}\nonumber
\end{align*}
where the last equality is due to condition \eqref{eq:a4}. \Cref{prop:moves_property2} specifies how many moves are made in this process for each type.
\begin{itemize}
\item Braid moves: 
\begin{align*}
\sum_{k=i+1}^{j-1}\left(1-M_{i,k+1}(w)+M_{i,k}(w)\right)=j-i-1-M_{i,j}(w)=j-i-1+\alpha_j-\alpha_i'.
\end{align*}
\item Commutation moves: 
\begin{align*}
\alpha_i'+\sum_{k=i+1}^{j-1}\left(\alpha_k -2(1-M_{i,k+1}(w)+M_{i,k}(w))\right)&=\sum_{k=i}^{j-1}\alpha_k' -2(j-i-1-M_{i,j}(w))\\
&=\sum_{k=i}^{j-1}\alpha_k' -2(j-i-1+\alpha_j-\alpha_i').
\end{align*}
\end{itemize}

Now, performing all the moves together, using condition \eqref{eq:a2}, 
\begin{align*}
\rrt{w}{I} &= \rr{w}_{1} \cdots \rr{w}_{i-1} \cdot (s_{\alpha_{i}+i-1} s_{\alpha_{i}+i-2} \cdots s_{\alpha_{i}'+i+1} )\cdot \rr{w}_{i}' \cdots \rr{w}_{j-1}' \cdot \rr{w}_{j}  \cdots \rr{w}_{n-1} \\
& = \rr{w}_{1} \cdots \rr{w}_{i-1} \cdot\rr{w}_{i}' \cdots \rr{w}_{j-1}' \cdot ( s_{\alpha_{j}'+j-1} s_{\alpha_{j}'+j-2} \cdots s_{\alpha_{j}+j} ) \cdot \rr{w}_{j}  \cdots \rr{w}_{n-1}= \rr{w'}.
\end{align*}

The process is applied $(\alpha_i-\alpha_i'-1)$ times and, therefore,
\begin{align*}
    \#\{k\colon f_k \mbox{ is a braid move}\} &= (\alpha_{i}-\alpha_{i}' - 1) (j-i-1+\alpha_{j}-\alpha_{i}'),\\
    \#\{k\colon f_k \mbox{ is a commutations}\} &= (\alpha_{i}-\alpha_{i}' - 1) \left(\left(\sum_{k=i}^{j-1} \alpha_{k}'\right) - 2 (j-i-1+\alpha_{j}-\alpha_{i}')\right).
\end{align*}

Now, suppose that $g_1 g_2 \cdots g_q$ is another sequence of moves from $\rrt{w}{I}$ to $\rr{w'}$. Then, the sequence $f_1 \cdots f_p g_q \cdots g_1$ is a loop in the graph of reduced decompositions of $w'$. By \Cref{prop:loopingw}, this loop has an even number of commutations and braid moves, which concludes the proof.
\end{proof}

\begin{remark} We claim that the construction in Proposition \ref{prop:countbraid} provides the shortest path in $R(w)$ between $\rrt{w}{I}$ into $\rr{w'}$.
\end{remark}

\subsection{Determination of degrees} 

The computation of the coefficient is based on the analysis of the effect of the two basic moves (commutation and braid) on equivalent reduced decompositions. 

\begin{lemma}\label{lem:degree_comm}
Let $\rr{w_{c}}$ be a reduced decomposition of $w$ such that $\rr{w_{c}}$ differs from $\rr{w}$ by a single commutation. Then, $\deg(\Phi_{\rr{w_{c}}}^{-1} \circ \Phi_{\rr{w}}) = -1$.
\end{lemma}
\begin{proof}

Denote $\ell = \ell(w)$, $\rr{w} = \cdots s_{i} s_{j} \cdots$, and $\rr{w_{c}} = \cdots s_{j} s_{i} \cdots$, where $|j-i| \geqslant 2$. More precisely, if $\rr{w} = r_{1} \cdots r_{\ell}$, where $r_{k} r_{k+1} = s_{i} s_{j}$, then $\rr{w_{c}} = r_{1}' \cdots r_{\ell}'$, where $r_{k}' r_{k+1}' = s_{j} s_{i}$ and $r_{l} = r_{l}'$ otherwise.

By \cite[Proposition 1.9]{RSm19}, the attaching maps  $\Phi_{\rr{w}}\colon [0,\pi]^\ell \rightarrow \mathcal{S}_w$ and $\Phi_{\rr{w_{c}}} \colon [0,\pi]^\ell \rightarrow \mathcal{S}_w$ are given by
\begin{align*}
\Phi_{\rr{w}}(t_{1},\dots, t_{\ell})&=g(t_{1},\dots, t_{k-1})\,e^{t_{k}A_{i}}\,e^{t_{k+1} A_{j}}\, h(t_{k+2},\dots, t_{\ell})\cdot b_0, \\
\Phi_{\rr{w_{c}}}(t_{1}',\dots, t_{\ell}')&=g(t_{1}',\dots, t_{k-1}')\,e^{t_{k}'A_{j}}\,e^{t_{k+1}' A_{i}}\, h(t_{k+2}',\dots, t_{\ell}')\cdot b_0,
\end{align*}
where $g(t_{1},\dots, t_{k-1}) = e^{t_{1}A_{r_{1}}}\cdots e^{t_{k-1}A_{r_{k-1}}}$ and $h(t_{k+2},\dots, t_{\ell}) = e^{t_{k+2}A_{r_{k+2}}}\cdots e^{t_{\ell}A_{r_{\ell}}}$.

Let us compute the degree of the map $\Phi^{-1}_{\rr{w_{c}}} \circ \Phi_{\rr{w}}$, which is a diffeomorphism $S^\ell \rightarrow S^\ell$ between the spheres obtained by collapsing the boundary of the cube $[0,\pi]^\ell$ to a single point (see \Cref{rem:degree}).
Since $|i-j|\geqslant 2$, the matrices $A_{i}$ and $A_{j}$ commute, that is, $e^{tA_{i}} e^{sA_{j}}= e^{sA_{j}} e^{tA_{i}}$, for $t,s \in [0,\pi]$. The map $\Phi_{\rr{w_{c}}}^{-1} \circ \Phi_\rr{w}$ restricted to the interior of $[0,\pi]^{\ell}$ is given by
\begin{equation*}
(\Phi_{\rr{w_{c}}}^{-1}\circ \Phi_{\rr{w}})(t_{1},\dots, t_{k}, t_{k+1}, \dots, t_{\ell}) = (t_{1},\dots, t_{k+1}, t_{k}, \dots, t_{\ell}).
\end{equation*}

The degree of $\Phi_{\rr{w_{c}}}^{-1}\circ \Phi_{\rr{w}}$ is the determinant of change of base map 
\begin{equation*}
(e_{1}, \dots, e_{k}, e _{k+1}, \dots, e_{\ell}) \mapsto (e_{1}, \dots, e _{k+1}, e_{k}, \dots, e_{\ell})
\end{equation*}
which is equal to $-1$.
\end{proof}

\begin{lemma}\label{lem:degree_braid}
Let $\rr{w_{b}}$ be a reduced decomposition of $w$ such that $\rr{w_{b}}$ differs from $\rr{w}$ by a single braid move. Then, $\deg(\Phi_{\rr{w_{b}}}^{-1} \circ \Phi_{\rr{w}}) = 1$.
\end{lemma}
\begin{proof}
Denote $\ell = \ell(w)$, $\rr{w} = \cdots s_{i} s_{i+1} s_{i} \cdots$, and $\rr{w_{b}} = \cdots s_{i+1} s_{i} s_{i+1} \cdots$. More precisely, if $\rr{w} = r_{1} \cdots r_{\ell}$, where $r_{k} r_{k+1} r_{k+2} = s_{i} s_{i+1} s_{i}$, then $\rr{w_{b}} = r_{1}' \cdots r_{\ell}'$, where $r_{k}' r_{k+1}' r_{k+2}' = s_{i+1} s_{i} s_{i+1}$ and $r_{l} = r_{l}'$ otherwise.

By \cite[Proposition 1.9]{RSm19}, the attaching maps  $\Phi_{\rr{w}}\colon [0,\pi]^\ell \rightarrow \mathcal{S}_w$ and $\Phi_{\rr{w_{b}}} \colon [0,\pi]^\ell \rightarrow \mathcal{S}_w$ are given by
\begin{align*}
\Phi_{\rr{w}}(t_{1},\dots, t_{\ell})&=g(t_{1},\dots, t_{k-1})\,e^{t_{k}A_{i}}\,e^{t_{k+1} A_{i+1}}\,e^{t_{k+2}A_{i}}\, h(t_{k+3},\dots, t_{\ell})\cdot b_0, \\
\Phi_{\rr{w_{b}}}(t_{1}',\dots, t_{\ell}')&=g(t_{1}',\dots, t_{k-1}')\,e^{t_{k}'A_{i+1}}\,e^{t_{k+1}' A_{i}}\,e^{t_{k+2}'A_{i+1}}\, h(t_{k+3}',\dots, t_{\ell}')\cdot b_0,
\end{align*}
where $g(t_{1},\dots, t_{k-1}) = e^{t_{1}A_{r_{1}}}\cdots e^{t_{k-1}A_{r_{k-1}}}$ and $h(t_{k+3},\dots, t_{\ell}) = e^{t_{k+3}A_{r_{k+3}}}\cdots e^{t_{\ell}A_{r_{\ell}}}$.

Let us compute the degree of the map $\Phi^{-1}_{\rr{w_{b}}} \circ \Phi_{\rr{w}}$, which is a diffeomorphism $S^\ell \rightarrow S^\ell$ between the spheres obtained by collapsing the boundary of the cube $[0,\pi]^\ell$ to a single point (see \Cref{rem:degree}). Since \cite[Proposition 1.9]{RSm19} says that the map $\Phi_{\rr{w}}$ is a diffeomorphism restricted to its interior, this will be done by checking the diffeomorphism at the coordinates curves in the interior of the cube. Our strategy here is to look at what happens to each curve of the form $(\pi/2, \ldots, \pi/2, t, \pi/2, \ldots, \pi/2), t \in (0,\pi)$.

For those coordinate curves where the reduced decompositions of $\rr{w}$ and $\rr{w_{b}}$ coincide, it follows directly that, for any integer $j \in [\ell] - \{k,k+1,k+2\}$
\begin{equation*}
\Phi_{\rr{w}}\left(\frac{\pi}{2},\dots, \frac{\pi}{2}, t_{j}, \frac{\pi}{2}, \cdots, \frac{\pi}{2}\right)=\Phi_{\rr{w_{b}}}\left(\frac{\pi}{2},\dots, \frac{\pi}{2}, t_{j}, \frac{\pi}{2}, \cdots, \frac{\pi}{2}\right) \, , \, t_{j} \in [0,\pi].
\end{equation*}

With respect to the remaining curves, by matrix calculus, the following equations for $t\in[0,\pi]$ are valid:
\begin{align}
e^{tA_{i}} e^{\frac{\pi}{2}A_{i+1}} e^{\frac{\pi}{2}A_{i}} &= e^{\frac{\pi}{2}A_{i+1}} e^{\frac{\pi}{2}A_{i}} e^{tA_{i+1}};\label{eq:tbraid1}\\
e^{\frac{\pi}{2}A_{i}} e^{tA_{i+1}} e^{\frac{\pi}{2}A_{i}} &= e^{\frac{\pi}{2}A_{i+1}} e^{(\pi-t)A_{i}} e^{\frac{\pi}{2}A_{i+1}};\label{eq:tbraid2}\\
e^{\frac{\pi}{2}A_{i}} e^{\frac{\pi}{2}A_{i+1}} e^{tA_{i}} &= e^{tA_{i+1}} e^{\frac{\pi}{2}A_{i}} e^{\frac{\pi}{2}A_{i+1}}.\label{eq:tbraid3}
\end{align}

Then, it follows respectively from Equations \eqref{eq:tbraid1},\eqref{eq:tbraid2} and \eqref{eq:tbraid3} that
\begin{align*}
\Phi_{\rr{w}}\left(\frac{\pi}{2},\dots, t_{k}, \frac{\pi}{2}, \frac{\pi}{2}, \cdots, \frac{\pi}{2}\right)&=\Phi_{\rr{w_{b}}}\left(\frac{\pi}{2},\dots, \frac{\pi}{2}, \frac{\pi}{2}, t_{k}, \cdots, \frac{\pi}{2}\right) \, , \, t_{k} \in [0,\pi];\\
\Phi_{\rr{w}}\left(\frac{\pi}{2},\dots, \frac{\pi}{2}, t_{k+1}, \frac{\pi}{2}, \cdots, \frac{\pi}{2}\right)&=\Phi_{\rr{w_{b}}}\left(\frac{\pi}{2},\dots, \frac{\pi}{2},\pi-t_{k+1},\frac{\pi}{2}, \cdots, \frac{\pi}{2}\right) \, , \, t_{k+1} \in [0,\pi];\\
\Phi_{\rr{w}}\left(\frac{\pi}{2},\dots, \frac{\pi}{2}, \frac{\pi}{2}, t_{k+2}, \cdots, \frac{\pi}{2}\right)&=\Phi_{\rr{w_{b}}}\left(\frac{\pi}{2},\dots, t_{k+2}, \frac{\pi}{2}, \frac{\pi}{2}, \cdots, \frac{\pi}{2}\right) \, , \, t_{k+2} \in [0,\pi].
\end{align*}

Therefore, the map $\Phi_{\rr{w_{b}}}^{-1}\circ \Phi_{\rr{w}}$ restricted to the interior of $[0,\pi]^{\ell}$ is
\begin{equation*}
(\Phi_{\rr{w_{b}}}^{-1}\circ \Phi_{\rr{w}})\left(t_{1},\dots, t_{k}, t_{k+1}, t_{k+2}, \dots, t_{\ell}) = (t_{1},\dots, t_{k+2},\pi-t_{k+1}, t_{k}, \dots, t_{\ell}\right).
\end{equation*}

The degree of $\Phi_{\rr{w_{b}}}^{-1}\circ \Phi_{\rr{w}}$ is the determinant of change of base map 
\begin{equation*}
(e_{1}, \dots, e_{k}, e _{k+1}, e_{k+2}, \dots, e_{\ell}) \mapsto (e_{1}, \dots, e_{k+2}, -e _{k+1}, e_{k}, \dots, e_{\ell})
\end{equation*}
which is equal to $1$.
\end{proof}

\subsection{The main formula}

With the details of the degree computations resolved, the scenario is now set to state and prove a definitive formula for computing the coefficients of the boundary map in the cellular homology of real flags of split real forms of type A. This formula eliminates any ambiguity in the choice of signs in the chain complex.

\begin{theorem}\label{thm:coeffA}
Let $w,w'\in S_{n}$ be such that $w'$ is covered by $w$, i.e., $w = w' \cdot (i,j)$, $1\leq i<j\leq n$. Assume that $\alpha = \code(w)$, $\alpha'= \code(w')$. Then, the number $I$ and the degree in \Cref{thm:rabelosanmartin} are given by
\begin{enumerate}
    \item [(i)] $I=\left(\sum_{k=1}^{i}\alpha_{k}\right)-\alpha_{i}'$;
    \item [(ii)] $\displaystyle \deg(\Phi_{\rrt{w}{I}}^{-1} \circ \Phi_{\rr{w'}}) = (-1)^{(\alpha_{i}-\alpha_{i}'-1)\cdot \sum_{l = i}^{j-1} \alpha_{l}'}$.
\end{enumerate}
\end{theorem}

\begin{proof}
The proof of \Cref{prop:countbraid} says that $\rrt{w}{I}$ comes from the deletion of the simple reflection $s_{\alpha_i'+i}$ from $\rr{w}'$ and
\begin{equation}
    \rrt{w}{I} = \rr{w}_{1} \cdots \rr{w}_{i-1} \cdot ( s_{\alpha_{i}+i-1} s_{\alpha_{i}+i-2} \cdots s_{\alpha_{i}'+i+1} \cdot \widehat{s}_{\alpha_{i}'+i} )\cdot \rr{w}_{i}' \cdot \rr{w}_{i+1}  \cdots \rr{w}_{n-1}.
\end{equation}

This occurs at position $I = \sum_{k=1}^{i-1} \alpha_k + (\alpha_i+i-1-(\alpha_i'+i))+1 = \sum_{k=1}^{i} \alpha_k -\alpha_i'$.

Let us prove the second assertion. Both reduced decompositions $\rrt{w}{I}$ and $\rr{w'}$ of $w'$ can be transformed into each other by a sequence of commutations and braid moves as described in the proof of \Cref{prop:countbraid}, i.e., there is a sequence of moves $f_1 f_2 \cdots f_p$ such that
\begin{equation*}
\rrt{w}{I} = \rr{w}^{0} \xrightarrow{f_1} \rr{w}^{1} \xrightarrow{f_2} \rr{w}^{2} \xrightarrow{f_3} \cdots \xrightarrow{f_p} \rr{w}^{r} = \rr{w'}.
\end{equation*}

If $\rr{w}^{i} $ transforms into $ \rr{w}^{i+1}$ by a single commutation, then Lemma \ref{lem:degree_comm} says that $ \deg(\Phi_{\rr{w}^{i}}^{-1} \circ \Phi_{\rr{w}^{i+1}}) = -1$. On the other hand, if $\rr{w}^{i} $ transforms into $ \rr{w}^{i+1}$ by a braid move, then Lemma \ref{lem:degree_braid} says that $ \deg(\Phi_{\rr{w}^{i}}^{-1} \circ \Phi_{\rr{w}^{i+1}}) = 1$.

Since $\deg(\Phi_{\rrt{w}{I}}^{-1} \circ \Phi_{\rr{w'}}) = \deg(\Phi_{\rr{w}^{0}}^{-1} \circ \Phi_{\rr{w}^{1}}) \cdot \deg(\Phi_{\rr{w}^{1}}^{-1} \circ \Phi_{\rr{w}^{2}}) \cdots \deg(\Phi_{\rr{w}^{r-1}}^{-1} \circ \Phi_{\rr{w}^{r}})$, then
\begin{equation*}
\deg(\Phi_{\rrt{w}{I}}^{-1} \circ \Phi_{\rr{w'}}) = (-1)^{\#\{k\colon f_k \mbox{ is a commutation}\}} \cdot 1^{\#\{k\colon f_k \mbox{ is a braid move}\}}.
\end{equation*}

Therefore, by Proposition \ref{prop:countbraid},
\begin{equation*}
\deg(\Phi_{\rrt{w}{I}}^{-1} \circ \Phi_{\rr{w'}}) = (-1)^{(\alpha_{i}-\alpha_{i}' - 1) \left(\left(\sum_{k=i}^{j-1} \alpha_{k}'\right) - 2 (j-i-1+\alpha_{j}-\alpha_{i}')\right)} = (-1)^{(\alpha_{i}-\alpha_{i}' - 1)\left(\sum_{k=i}^{j-1} \alpha_{k}'\right)}.
\end{equation*}

Notice that if we choose any other sequence of moves from $\rrt{w}{I}$ to $\rr{w'}$ then \Cref{prop:countbraid} guarantees that the degree will not change since the number of commutations will have the same parity.
\end{proof}
 
The results in \Cref{thm:coeffA} have an insightful interpretation the code diagrams of $w$ and $w'$. Recall that the code diagram of $w'$ is obtained from the code diagram of $w$ by removing one box in the $i$-th row. 

Let $B_0,B_1$ and $B_2$ be given as
\begin{itemize}
\item $B_0=\sum_{k=1}^{i-1}\alpha_{k}$: the number of boxes below the $i$-th row;
\item $B_1=\alpha_i-\alpha'_i-1$: the number of boxes in the $i$-th row at right of the removed box;
\item $B_2=\sum_{k=i}^{j-1}\alpha'_{k}$: the number of boxes between the removed box (considering eventually the boxes at left) and the $j$-th row.
\end{itemize}

\begin{corollary}
If $c(w,w')\neq 0$, then $c(w,w')=(-1)^{(B_0+B_1+1)+B_1\cdot B_2 }\cdot 2$.
\end{corollary}

\begin{corollary}\label{coro:coef}
    $c(w,w')=(-1)^{(B_0+B_1+1)+B_1\cdot B_2 }\cdot (1+(-1)^{j-i})$.
\end{corollary}

\begin{figure}[hbtp]
    \centering
    \includegraphics[scale=0.7]{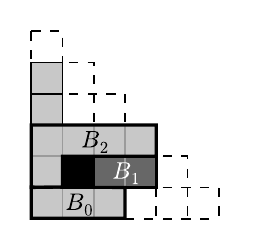}
    \caption{The coefficient computation using data from the diagram}
    \label{fig:ex_coeficiente}
\end{figure}

As before, let $w=4\, 6\, 7\, 2\, 3\, 1\, 5$ and $w'=4\, 2\, 7\, 6\, 3\, 1\, 5$ with the corresponding codes $\alpha = (3,4,4,1,1,0)$ and $\alpha' = (3,1,4,3,1,0)$. It follows that $w$ covers $w'$ with $w = w'\cdot (2,4)$. This example is illustrated in \Cref{fig:ex_sign} where we have shown how the black box is removed from $w$ to get $w'$.
By \Cref{thm:rabelosanmartin}, $c(w,w')=\pm 2$ since the difference $j-i=4-2=2$ is even. 
Now, observe in \Cref{fig:ex_coeficiente} that $B_0=3$, $B_1=2$ and $B_2=5$. Hence, by \Cref{coro:coef}
\begin{equation*}
c(w,w')=(-1)^{(3+2+1)+2\cdot5}\cdot(1-(-1)^2)=(-1)^{16}2=+2.
\end{equation*}  

\begin{figure}[ht]
\centering
\includegraphics[scale=0.8]{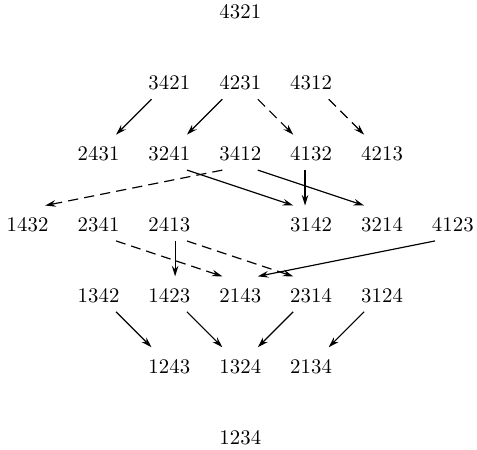}
\caption{Signed incidence diagram for the maximal flag manifold of type $A_{3}$.}
\label{fig:n4incidence}
\end{figure}

According to \cite{Mat19}, the incidence graph has Schubert cells $w\in \weyl$ as its vertices, and two vertices $w,w'$ are connected if $c(w,w')=\pm 2$. The signed incidence diagram provides additional information on the signs: a dashed line represents $+2$, while a solid line represents $-2$. For instance, using \Cref{thm:coeffA}, the coefficients for the maximal flag manifold $\mathbb{F}(1,2,3,4)$ are shown in Figure \ref{fig:n4incidence}. It may be set side by side with Figure 3 in \cite{Mat19}, where a different formula is used for the signs. Both formulas depend on a complex combinatorics of the permutations, but the formula in \Cref{thm:coeffA} has the advantage of being interpretable through the corresponding code diagrams.

\section{Low-dimensional homology groups for flag manifolds of type A}\label{sec:generators}

In the sequence, all generators of integral homology groups up to dimension $4$ and free part generators up to dimension 6 are explicitly computed using the techniques presented in this paper. From now on, denote the Schubert cell $\schub_{w}$ only by its index $w$ to avoid repetition and simplify the notation.

\subsection{Code Spectrum}
Following \cite{LR22a}, given $w\in S_{n}$, the \emph{code spectrum} of $w$ is the unique partition $0<b_{1}\leqslant b_{2}\leqslant \cdots \leqslant b_{l}<n$ such that the code $\alpha$ of $w$ is given by $\alpha_{i} = \#\{j \colon b_{j} = i\}$.  Denote by $\w{b_{1},\cdots, b_{\ell}}$ the permutation $w$ given by this code spectrum to distinguish it from the other notations. The code spectrum is easily obtained from the diagram since it records the row of each box. For example, $w=1\,3\, 7 \,5\,8\,2\,9\,4\,6\in \weyl^{\Theta}$ corresponds to $\w{2,3,3,3,3,4,4,5,5,5,7,7}$.

This notation simplifies the description of permutations with a specific length. The following lemmas characterize the elements of the $2$- and $3$- skeleton parametrized by $\weyl^{\Theta}$ with respect to the code spectrum.

\begin{lemma}\label{lem:wthetalow1and2}[\cite{LR22a}, Lemma 5.6]
\begin{enumerate}
\item Given $i\in [n-1]$, $\w{i}\in \weyl^{\Theta}$ if, and only if, $a_{i}\not\in\Theta$;\
\item Given $i\in[n-2]$ and $j\in [n-1]$ such that $i\leqslant j$, $\w{i,j} \in \weyl^{\Theta}$ if, and only if, one of the following happens:
\begin{itemize}
\item $a_{i},a_{j}\not\in\Theta$;
\item $j=i+1$, $a_{i}\in \Theta$, and $a_{i+1}\not\in\Theta$.\
\end{itemize}
\end{enumerate}
\end{lemma}

\begin{lemma} Given $i \in [n-3], j \in [n-2]$ and $k\in [n-1]$ such that $i\leq j\leq k$, $\w{i,j,k}\in \weyl^{\Theta}$ if, and only if, one of the following happens:
\begin{itemize}
\item $a_i,a_j,a_k \not\in\Theta$;
\item $j=i+1$, $a_i \in \Theta$ and $a_{i+1},a_k \not\in\Theta$;
\item $k=j+1$, $a_j \in \Theta$ and $a_i, a_{j+1} \not\in\Theta$;
\item $k-1=j=i+1$, $a_i,a_{i+1} \in \Theta$ and $a_{i+2} \not\in\Theta$.
\end{itemize}
\end{lemma}
\begin{proof} The result follows as an application of the Proposition \ref{prop:bijection_wk}. For each choice of $\Theta$, arrange the boxes of the $3$-cells such that they form a partition inside the corresponding rectangles.

The possible types of $3$-cells are given by
\begin{align*}
s_{i}s_{j} s_{k} &= \w{i,j,k} \mbox{, for } i,j,k\in [n-1];\\
s_{i+1}s_{i} s_{k} &= \w{i,i,k} \mbox{, for } i\in [n-2],k\in [n-1];\\
s_{i}s_{j+1} s_{j} &= \w{i,j,j} \mbox{, for } i,j\in [n-2];\\
s_{i+2}s_{i+1} s_{i} &= \w{i,i,i} \mbox{, for } i \in [n-3].
\end{align*}

The possible configurations in which the $3$-cells form a partition are shown in Figure \ref{fig:3cells1case}.
\begin{figure}[hbtp]
\centering
\includegraphics[scale=0.65]{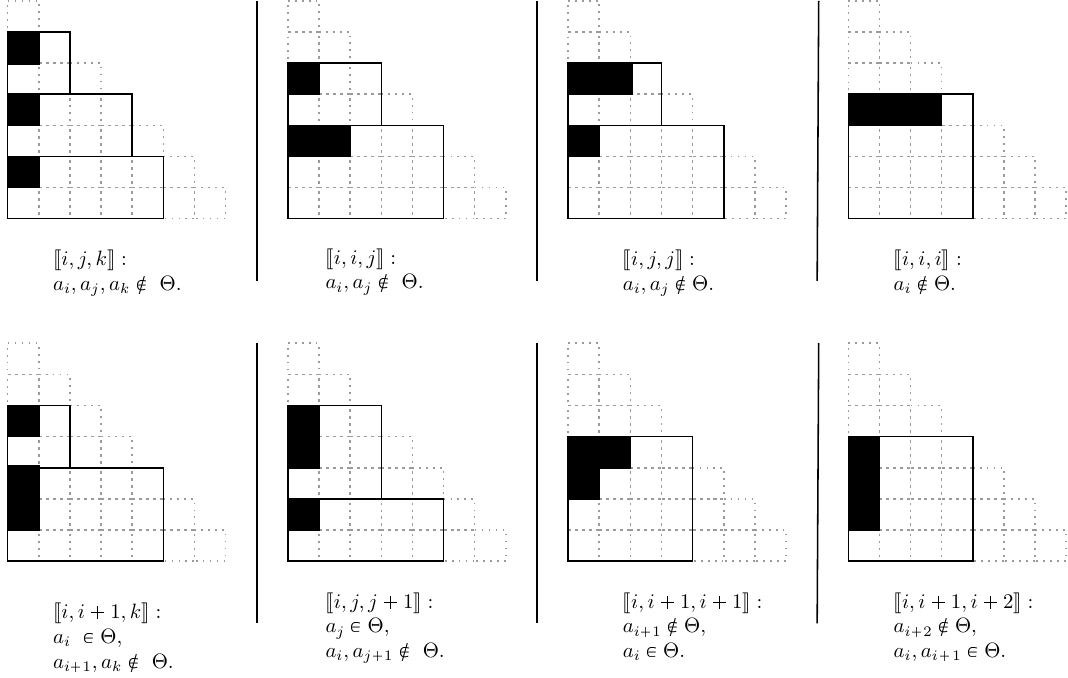}
\caption{$3$-cells diagrams. The vertical gap between the boxes may not exist since $i\leqslant j \leqslant k$}
\label{fig:3cells1case}
\end{figure}
\end{proof}

\begin{remark} It is not difficult to derive analogues of these lemmas for the $n$-skeleton with $n\geq 4$. Some of them will be used in the sequence of the paper. 
\end{remark}

\subsection{Low dimensional free generators}

Denote the number of connected components of the Dynkin diagram of $\Theta$ by $r=r(\Theta)$ and the number of connected components of the Dynkin diagram of $\Theta$ with exactly $k$ elements by $r_{k}= r_{k}(\Theta)$, for $k>0$. Also, denote $r_0=1$ if $\Theta$ is empty and $r_0=0$ otherwise.

The $1$- and $2$-homology groups are given by the next theorem.

\begin{theorem}[\cite{LR22a}, Theorem 5.7]\label{thm:homology1-2} Consider a partial flag manifold $\mathbb{F}_{\Theta}$ of $G=\mathrm{Sl}(n,\mathbb{R})$, where $\Theta\subset\Sigma= \{a_1, \ldots, a_{n-1}\}$.
\begin{enumerate}
\item For $n\geqslant 2$, the first homology group is
\begin{equation*}
H_{1}(\mathbb{F}_{\Theta},\Z) \cong (\Z_{2})^{n-|\Theta|-1}
\end{equation*}
and it is generated by the set of Schubert cells $\schub_{\w{i}}$ such that $a_{i}\in \Sigma-\Theta$.

\item For $n\geqslant 3$, the second homology group is
\begin{equation*}
H_{2}(\mathbb{F}_{\Theta},\Z) \cong (\Z_{2})^{T_2}
\end{equation*}
where $T_2 = \binom{n-|\Theta|-1}{2} +r - 1$ and it is generated by the set of Schubert cells
\begin{itemize}
\item $X_{i,j}=\schub_{\w{i,j}}$, for any $i,j\in [n-1]$ such that $j-i\geqslant 2$ and $a_{i},a_{j}\in \Sigma-\Theta$;
\item $X_{i,i+1}=\schub_{\w{i,i+1}} - \schub_{\w{i+1,i+1}}$, for every $i\in [n-3]$ such that $a_{i+1}\in \Sigma-\Theta$.
\end{itemize}
\end{enumerate}
\end{theorem}

In \cite{LR22a}, there is no mention that the formulas of $H_1$ and $H_2$ also work, respectively, for $n=2$ and $n=3$.

Define the indicator function of $i\in [n-1]$ with respect to $\Theta$ as follows
\begin{align*}
\mathds{1}^{\Theta}_i = \left\{ \begin{array}{ll} 1, & a_i \in \Theta; \\ 0, & a_i \not\in\Theta. \end{array}\right.
\end{align*}

Recall that $\notInTheta=\Sigma-\Theta$. It follows that $\mathds{1}^{\notInTheta}_i = 1 - \mathds{1}^{\Theta}_i$.

The following construction will be required to get the 4-homology group. Given any $p,q\in [n]$, with $p<q$, consider the sum
\begin{equation}
Z(p,q) = \sum_{i = p}^{q-1} {(\mathds{1}^{\notInTheta}_{i+1}) \w{i,i,i+1,i+1}}.
\end{equation}

Suppose that the number of connected components of $\Theta$ is at least two, that is, $r\geqslant 2$. Denote each component of $\Theta$ by $\Theta_k =\{a_{t_k(1)},\dots, a_{t_k(u(k))}\}$ following the usual order of the roots, i.e.,
if $\Theta_k$ has $u(k)$ roots, then $t_i(j)$ is the index of the $j$-th root inside the $i$-th connected component and
\begin{equation*}
    t_1(1)<\cdots<t_1(u(1))<\cdots < t_k(1)<\cdots<t_k(u(k))<\cdots < t_{r}(1)<\cdots<t_{r}(u(r)).
\end{equation*}

For each $i\in [r-1]$, define
\begin{equation*}
Z_{i} = Z(t_i(u(i)), t_{i+1}(1)-1).
\end{equation*}

Let us illustrate $Z_i$ using the case where $n=8$ and $\Theta = \{a_1,a_2,a_6,a_8\}$. This choice is represented in the corresponding Dynkin diagram as

\vspace{10pt}

\begin{center}
    \begin{tikzpicture}[scale=2]
        \dynkin[labels={a_1,a_2,,,,a_6,,a_8}] A{ooxxxoxo}
    \end{tikzpicture}
\end{center}

\begin{figure}[ht]
    \centering
    \includegraphics[scale=0.6]{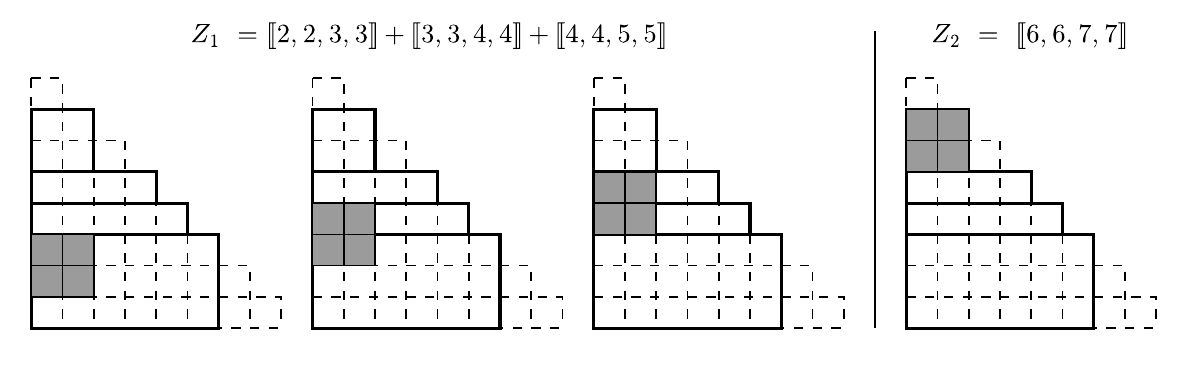}
    \caption{$Z_1$ and $Z_2$ for $n=8$ and $\Theta = \{a_1,a_2,a_6,a_8\}$.}
    \label{fig:h4ex}
\end{figure}

Hence $\Theta$ has three connected components, i.e., $r=r(\Theta)=3$, so that $Z_1=Z(2,5)=\w{2,2,3,3}+\w{3,3,4,4}+\w{4,4,5,5}$ and $Z_2=Z(6,7)=\w{6,6,7,7}$ as in \Cref{fig:h4ex}. With respect to the Dynkin diagram, each space between the connected components of $\Theta$ provides the summation of $Z_1$ and $Z_2$.

According to the next proposition, the Betti number $\beta_4=2$, i.e., these are the two generators for the $4$-homology of $\mathbb{F}(3,4,5,7,9)$.
\begin{proposition}\label{prop:betti_generators}
    Let $\Theta\subset\Sigma$. The nonzero Betti numbers up to $\beta_6$ with the corresponding generators in $H_i(\mathbb{F}_\Theta,\Z)$ are given by \Cref{tbl:homology_free}.

\begin{table}[ht]
    \centering
    \caption{Nonzero Betti numbers up to $\beta_6$ and the corresponding generators in $H_i(\mathbb{F}_\Theta,\Z)$}\label{tbl:homology_free}
    \begin{tabular}{llc}
        \toprule
        & Conditions &  Generators \\
        \midrule
        $\beta_3=2$ & $n = 4$ and $\Theta=\emptyset$ & $\w{1,1,2}; \w{1,1,1} + \w{1,2,3}$ \\
        $\beta_3=1$ & $n = 4$ and $r=1$ & $\one^{\notInTheta}_{1}\w{1,1,1} + \one^{\notInTheta}_{3}\w{1,2,3}$ \\
        $\beta_3=1$ & $n\neq 4$ and $\Theta=\emptyset$ & $\w{1,1,2}$ \\
        $\beta_4=r-1$ & $n\geqslant 4$ & $Z_i$, for $i\in[r-1]$ \\
        $\beta_5=1$ & $n = 6$ and $L\leqslant 2$ & $Y_5=\one^{\notInTheta}_{1}\w{1,1,1,1,1}+\one^{\notInTheta}_{3}\w{1,2,3,3,3}+$ \\
        & & $+\one^{\notInTheta}_{1}\one^{\notInTheta}_{5}\w{1,1,1,4,5}+\one^{\notInTheta}_{5}\w{1,2,3,4,5}$\\
        $\beta_6=1$ & $n=4$ and $\Theta=\emptyset$ & $\w{1,1,1,2,2,3}$ \\ 
        \bottomrule
    \end{tabular}
\end{table}
\end{proposition}
\begin{proof}
The proof consists of two steps: expanding the Poincaré polynomial to obtain the Betti numbers using the formulas given in \Cref{prop:poincare}, and verifying that the given cells generate the free part.

The cases for $\beta_3$ will be proved in Theorem \ref{thm:H3}.

If $n=4$ and $\Theta=\emptyset$ then $\beta_6=1$ and it corresponds to the principal involution in Figure \ref{fig:n4incidence}.

If $n=6$ and $L\leqslant 2$, then $\beta_5=1$. Notice that
\begin{align*}
    \one^{\notInTheta}_{1} \partial \w{1,1,1,1,1} &= -2(\one^{\notInTheta}_{1}\one^{\notInTheta}_{3})\w{1,3,3,3} -2(\one^{\notInTheta}_{1}\one^{\notInTheta}_{5})\w{1,1,1,5}; \\
    \one^{\notInTheta}_{3}\partial \w{1,2,3,3,3} &= 2(\one^{\notInTheta}_{1}\one^{\notInTheta}_{3})\w{1,3,3,3} -2(\one^{\notInTheta}_{3}\one^{\notInTheta}_{5})\w{1,2,3,5}; \\
    \one^{\notInTheta}_{1}\one^{\notInTheta}_{5}\partial \w{1,1,1,4,5} &= -2(\one^{\notInTheta}_{1}\one^{\notInTheta}_{5})\w{1,3,4,5} +2(\one^{\notInTheta}_{1}\one^{\notInTheta}_{5})\w{1,1,1,5}; \\
    \one^{\notInTheta}_{5}\partial \w{1,2,3,4,5} &= 2(\one^{\notInTheta}_{1}\one^{\notInTheta}_{5})\w{1,3,4,5} +2(\one^{\notInTheta}_{3} \one^{\notInTheta}_{5})\w{1,2,3,5}. 
\end{align*}

Hence, for 
\begin{align*}
Y_5=\one^{\notInTheta}_{1}\w{1,1,1,1,1}+\one^{\notInTheta}_{3}\w{1,2,3,3,3}+\one^{\notInTheta}_{1}\one^{\notInTheta}_{5}\w{1,1,1,4,5}+\one^{\notInTheta}_{5}\w{1,2,3,4,5},
\end{align*}
$\partial Y_5=0$. Moreover, if $w\in \{\w{1,1,1,1,1}, \w{1,2,3,3,3}, \w{1,1,1,4,5}, \w{1,1,1,4,5}\}$ and $v$ covers $w$ then  $c(v,w)=0$. Therefore, $Y_5$ generates the free part of $H_5(\mathbb{F}_\Theta,\mathbb{Z})$.

Let us carefully study the 4-homology. Notice that
\begin{align*}
\partial Z(p,q) &= (\mathds{1}^{\notInTheta}_{p+1}) (2 (\mathds{1}^{\notInTheta}_{p+2}){\w{p+1,p+1,p+2}}-2 (\mathds{1}^{\notInTheta}_{p}){\w{p,p,p+1}})) + \\
&+ (\mathds{1}^{\notInTheta}_{p+2}) (2 (\mathds{1}^{\notInTheta}_{p+3}){\w{p+2,p+2,p+3}}-2 (\mathds{1}^{\notInTheta}_{p+1}){\w{p+1,p+1,p+2}}))+\\
&+\dots+\\
&+ (\mathds{1}^{\notInTheta}_{q-1}) (2(\mathds{1}^{\notInTheta}_{q}){\w{q-1,q-1,q}}-2 (\mathds{1}^{\notInTheta}_{q-2}){\w{q-2,q-2,q-1}}))+\\
&+ (\mathds{1}^{\notInTheta}_{q}) (2(\mathds{1}^{\notInTheta}_{q+1}){\w{q,q,q+1}}-2 (\mathds{1}^{\notInTheta}_{q-1}){\w{q-1,q-1,q}})).
\end{align*}

Suppose that the number of connected components of $\Theta$ is at least two, that is, $r\geqslant 2$. Clearly, $Z_i\neq 0$ and $\partial Z_i=0$ for $i\in [r-1]$. Moreover, they are linearly independent.

If $w\in\weyl^{\Theta}$ covers $\w{i,i,i+1,i+1}$ then there is $j\in [n-1]$ such that either $w=\w{j,i,i,i+1,i+1}$ or $w=\w{i,i,i,i+1,i+1,j}$ and therefore $c(w, \w{i,i,i+1,i+1})=0$.

Let us check that $\{Z_i\colon i\in [r-1]\}$ generates the free part. According to \Cref{prop:poincare}, in the formula for the Poincaré polynomial, the Betti number $\beta_4$ is provided exclusively by the $t$-multinomial part. If $r \leqslant 1$ then $L = 0$ and the $t$-multinomial is trivial, with the result that $\beta_4=0$.

Assume $r \geqslant 2$. For $\Theta'= \{a_{k_{1}},\dots, a_{k_{s}}\}$, if $\left\lfloor \frac{k_{j}-k_{j-1}}{2} \right\rfloor = 0$, then $k_{j-1}$ and $k_{j}$ are consecutive integers. Let $\alpha_1,\dots, \alpha_u$ be the non-zero numbers of the list $\left\lfloor \frac{k_{1}}{2} \right\rfloor , \left\lfloor \frac{k_{2}-k_{1}}{2} \right\rfloor , \cdots , \left\lfloor \frac{n-k_{s}}{2} \right\rfloor$. In this case, $u$ is the number of connected components of $\Theta$, i.e., $u=r$.

Let $M=\{1^{\alpha_{1}}, 2^{\alpha_{2}},\dots, r^{\alpha_{r}}\}$ be the multiset of cardinality $L$. By \cite[Proposition 1.7.1]{Sta11}, 
\begin{equation*}
    \binom{L}{ \left\lfloor \frac{k_{1}}{2} \right\rfloor , \left\lfloor \frac{k_{2}-k_{1}}{2} \right\rfloor , \cdots , \left\lfloor \frac{n-k_{s}}{2} \right\rfloor}_{t^{4}} = \sum_{w\in S_M} t^{4 \ell(w)}.
\end{equation*}

Then, $\beta_4 = \#\{w\in S_M\colon \ell(w)=1\}=r-1$.
\end{proof}

Notice that the above results match with those of Casian and Kodama in \cite{CK2010} about the generators for the cohomology ring of real Grassmannian manifolds. For instance, the Betti number $\beta_4=1$ for $\mathrm{Gr}(2,4)$ is generated by $Z_1=\w{1,1,2,2}$ as obtained in \cite{CK2010}, Example 5.1; the Betti number $\beta_5=1$ for $\mathrm{Gr}(3,6)$ is generated by $\w{1,2,3,3,3}$ as it is obtained in \cite{CK2010}, Example 5.2.


\subsection{Torsion of the third and fourth homology groups}

To compute the torsion, a first approach is given by \Cref{lem:chenhe}, which provides a formula for the torsion part in terms of the free part. 
This approach requires knowledge of the free part and would lead to a case-by-case analysis. Using the combinatorial methods developed in this text, explicit formulas for the torsion of the third and fourth homology groups are obtained in terms of $\Theta$ and its connected components, including the expression for the generators. It does not seem to be possible to obtain these results from the former approach.


\begin{theorem}\label{thm:H3} The dimension of torsion for 3-homology is
\begin{itemize}
    \item For $n=3$: $T_3 = 0$;
    \item For $n=4$: $T_3 = 2-|\Theta|$;
    \item For $n\geqslant 5$:
    \begin{equation*}
    T_{3} = \binom{n-|\Theta|}{3} + r (n - |\Theta|-1) - r_{0}  - r_{1}.
    \end{equation*}
\end{itemize}

All generators will be presented along the proof.
\end{theorem}
\begin{proof}
Suppose that $\w{i,j,k}\in\weyl^{\Theta}$. If $n=3$, then $\w{1,1,2}$ is the only element in the kernel when $|\Theta|=0$.

If $n=4$ then 
\begin{align*}
    \partial \w{1,1,1} &= -2 (\one_{3}^{\notInTheta}) \w{1,3}; &
    \partial \w{1,1,2} &= 0;\\
    \partial \w{1,1,3} &= 0;&
    \partial \w{1,2,2} &= 2 \w{1,2} - 2 \w{2,2};\\
    \partial \w{1,2,3} &= 2 (\one_{1}^{\notInTheta}) \w{1,3};&
    \partial \w{2,2,3} &= 0.
\end{align*}

The kernel of $\partial$ is generated by $X_{1} = \one^{\notInTheta}_{1}\w{1,1,1} + \one^{\notInTheta}_{3}\w{1,2,3}, X_{1,1,2} = \w{1,1,2}, X_{1,1,3} = \w{1,1,3}$ and $X_{2,2,3} = \w{2,2,3}$. The generators of the homology are presented in \Cref{tbl:genH3_n4}.

\begin{table}[ht]
    \centering
    \caption{Generators of the 3-homology for $n=4$}\label{tbl:genH3_n4}
    \begin{tabular}{cccccccc}
        \toprule
        $\Theta$ & $\emptyset$ & $\{a_1\}$ & $\{a_2\}$ & $\{a_3\}$ & $\{a_1,a_2\}$ & $\{a_2,a_3\}$ & $\{a_1,a_3\}$\\ \midrule
        Free part & $X_1, X_{1,1,2}$ & $X_1$ & $X_1$ & $X_1$ & $X_1$ & $X_1$ & 0 \\ 
        Torsion & $X_{1,1,3}, X_{2,2,3}$ & $X_{2,2,3}$ & $X_{1,1,3}$ & $X_{1,1,2}$ & 0 & 0 & 0 \\ \bottomrule
    \end{tabular}
\end{table}

If $n\geqslant 5$ then the boundary map of every $3$-cell in the flag manifold $\mathbb{F}_{\Theta}$ is given by
\begin{align}
&\partial{\w{i,i,i}} = -2(\mathds{1}^{\notInTheta}_{i+2}){\w{i,i+2}}, \  i\in [n-3];\label{eq:b3-1}\\
&\textcolor{blue}{\partial{\w{i,i,i+1}} = 0, \  i\in [n-2];^\star}\label{eq:b3-2}\\
&\textcolor{blue}{\partial{\w{i,i,i+2}} = 0, \  i\in [n-3];^\star}\label{eq:b3-3}\\
&\partial{\w{i,i,k}} = -2{\w{i,k}}, \  i\in [n-4],\;\; k\in [i+3,n-1];\label{eq:b3-4}\\
&\textcolor{red}{\partial{\w{i,i+1,i+1}} = 2{\w{i,i+1}} - 2{\w{i+1,i+1}}, \;\; i\in [n-3];^\bullet}\label{eq:b3-5}\\
&\partial{\w{i,k,k}} = 2{\w{i,k}}, \  i\in [n-4] , \;\; k\in [i+2,n-2];\label{eq:b3-6}\\
&\partial{\w{i,i+1,i+2}} = 2(\mathds{1}^{\notInTheta}_{i}){\w{i,i+2}}, \  i\in [n-3];\label{eq:b3-7}\\
&\partial{\w{i,k-1,k}} = 2{\w{i,k}}, \  i\in [n-4] ,\;\; k\in [i+3,n-1];\label{eq:b3-8}\\
&\partial{\w{i,i+1,k}} = -2{\w{i+1,k}}, \  i\in [n-4] ,\;\; k\in [i+3,n-1];\label{eq:b3-9}\\
&\textcolor{blue}{\partial{\w{i,j,k}} = 0, \  i\in [n-5],\;\; j\in [i+2,n-3],\;\; k\in [j+2,n-1].^\star}\label{eq:b3-10}
\end{align}

Let us compute the generators of $\ker(\partial)$ that will be obtained by combining the boundary map in Equations \eqref{eq:b3-1} to \eqref{eq:b3-10}.
The formulas above suggest that there are three types of boundary maps (marked with different symbols).

Equations \eqref{eq:b3-2}, \eqref{eq:b3-3}, and \eqref{eq:b3-10} marked with $\star$ provide the generators:
\begin{align*}
& X_{1,1,2}=\w{1,1,2};\\
& X_{i,i,i+1}=\w{i,i,i+1} - \mathds{1}^{\notInTheta}_{i-1}\w{i-1,i-1,i}, \  i\in [2,n-2];\\
& X_{i,i,i+2}=\w{i,i,i+2}, \  i\in [n-3];\\
& X_{i,j,k}=\w{i,j,k}, \  i\in [n-5],\;\; j\in [i+2,n-3],\;\; k\in [j+2,n-1].
\end{align*}

By Equation \eqref{eq:b3-5} marked with $\bullet$, no combination of Schubert cells with $\w{i,i+1,i+1}$ belongs to the kernel of $\partial$ since it is not possible to eliminate the term $\w{i+1,i+1}$.

For the remaining Equations \eqref{eq:b3-1}, \eqref{eq:b3-4}, \eqref{eq:b3-6}, \eqref{eq:b3-7}, \eqref{eq:b3-8}, and \eqref{eq:b3-9}, proceed finding the options of $i$ and $k$ such that the boundaries can be combined with each other. \Cref{tbl:h3_1} summarizes such information.

\begin{table}[ht]
\centering
\caption{Comparing equations without marker}\label{tbl:h3_1}
\begin{tabular}{cccc}
\toprule
 & $k=i+2$ & $k\in [i+3,n-2]$ & $k=n-1$ \\ 
\midrule
$i=1$ & \eqref{eq:b3-1} \eqref{eq:b3-6} \eqref{eq:b3-7} & \eqref{eq:b3-4} \eqref{eq:b3-6} \eqref{eq:b3-8} & \eqref{eq:b3-4} \eqref{eq:b3-8} \\ 
$i\in [2, n-5]$ & \eqref{eq:b3-1} \eqref{eq:b3-6} \eqref{eq:b3-7} \eqref{eq:b3-9} & \eqref{eq:b3-4} \eqref{eq:b3-6} \eqref{eq:b3-8} \eqref{eq:b3-9} & \eqref{eq:b3-4} \eqref{eq:b3-8} \eqref{eq:b3-9} \\ 
$i=n-4$ & \eqref{eq:b3-1} \eqref{eq:b3-6} \eqref{eq:b3-7} \eqref{eq:b3-9} & \eqref{eq:b3-4} \eqref{eq:b3-6} \eqref{eq:b3-8} \eqref{eq:b3-9} & \eqref{eq:b3-4} \eqref{eq:b3-8} \eqref{eq:b3-9}\\ 
$i=n-3$ & \eqref{eq:b3-1} \eqref{eq:b3-7} \eqref{eq:b3-9} & \eqref{eq:b3-9} & \eqref{eq:b3-9} \\ 
\bottomrule 
\end{tabular} 
\end{table}

The generators are given by the following combinations:
\begin{align*}
\mbox{\eqref{eq:b3-1}}+\mbox{\eqref{eq:b3-7}: } & X_{i,i,i} = \w{i,i,i} + \mathds{1}^{\notInTheta}_{i+2} \w{i,i+1,i+2}, \  i\in[n-3];\\
\mbox{\eqref{eq:b3-7}}-\mbox{\eqref{eq:b3-6}: } & X_{i,i+1,i+2} = \w{i,i+1,i+2} - \mathds{1}^{\notInTheta}_{i}\w{i,i+2,i+2}, \  i\in[n-4]; \\
\mbox{\eqref{eq:b3-9}}+\mbox{\eqref{eq:b3-7}: } & X_{i-1,i,i+2} = \w{i-1,i,i+2} + \w{i,i+1,i+2}, \  i\in[2,n-3]\\
\mbox{or } & X_{i,i+1,i+3} = \w{i,i+1,i+3} + \w{i+1,i+2,i+3}, \  i\in[1,n-4];\\
\mbox{\eqref{eq:b3-4}}+\mbox{\eqref{eq:b3-8}: } & X_{i,i,k} = \w{i,i,k} + \w{i,k-1,k}, \  i\in[n-4], k\in [i+3,n-1];\\
\mbox{\eqref{eq:b3-8}}-\mbox{\eqref{eq:b3-6}: } & X_{i,k-1,k} = \w{i,k-1,k} - \w{i,k,k}, \  i\in [n-4], k\in[i+3,n-2];\\
\mbox{\eqref{eq:b3-9}}+\mbox{\eqref{eq:b3-8}: } & X_{i-1,i,k} = \w{i-1,i,k} + \w{i,k-1,k}, \  i\in [2,n-4], k \in [i+3,n-1]\\
\mbox{or } & X_{i,i+1,k} = \w{i,i+1,k} + \w{i+1,k-1,k}, \  i\in [n-5], k \in [i+4,n-1].
\end{align*}

If $\Theta = \emptyset$, then from \Cref{tbl:homology_free}, $X_{1,1,2}$ is the only generator of the free part. Then, all elements different from $X_{1,1,2}$ are the image through $\partial$ and generate the torsion.


If $\Theta \neq \emptyset$, then all the generators of $\ker(\partial)$ also generate the torsion since there is no free part.

The torsion generators will be counted by splitting up according to the following disjoint sets:
\begin{align*}
    A_1 &=\{X_{i,i,i} \colon a_i \not\in\Theta\}; & A_4 &=\{X_{i,i+1,k} \colon a_i \in \Theta; a_{i+1},a_k \not\in\Theta\};\\
    A_2 &=\{X_{i,i,k} \colon a_i,a_k \not\in\Theta\}; & A_5 &=\{X_{i,k-1,k} \colon a_{k-1} \in \Theta; a_{i},a_k \not\in\Theta\};\\
    A_3 &=\{X_{i,j,k} \colon a_i,a_k,a_j \not\in\Theta\}; & A_6 &=\{X_{i,i+1,i+2} \colon a_i,a_{i+1} \in \Theta; a_{i+2} \not\in\Theta\}.
\end{align*}

It follows that:
\begin{itemize}
\item For $|A_{1}|$: Since $X_{i,i,i}$ is defined for $i\in [n-3]$, exclude the choices for $a_i \in \Sigma - \Theta$ in the interval $[n-3]$. Then,
\begin{align*}
|A_{1}| &= \left\{
\begin{array}{cl}
n-3 -|\Theta|, & \mbox{if } a_{n-2},a_{n-1} \not\in\Theta; \\ 
n-3 -(|\Theta|-1), & \mbox{if } a_{n-2}\in \Theta \mbox{ and } a_{n-1} \not\in\Theta; \\ 
n-3 -(|\Theta|-1), & \mbox{if } a_{n-2}\not\in\Theta \mbox{ and } a_{n-1} \in \Theta; \\ 
n-3 -(|\Theta|-2), & \mbox{if } a_{n-2}, a_{n-1} \in \Theta; \\ 
\end{array} 
\right.\\
&= n-3-|\Theta|+\mathds{1}_{n-2}^{\Theta}+\mathds{1}_{n-1}^{\Theta}.
\end{align*}

\item For $|A_{2}|$: Since $X_{i,i,k}$ are defined for $i\in [n-2]$ and $k\in[i+1,n-1]$, exclude $|\Theta|$ choices among all possible choices of pairs $i<k$ in the interval $[n-1]$. Then,
\begin{equation*}
|A_{2}| = \binom{n-1-|\Theta|}{2}.
\end{equation*}

\item For $|A_{3}|$: There are $n-3-|\Theta|$ choices for $i$ such that $a_i \not\in\Theta$. However, if $a_{n-2},a_{n-1} \in \notInTheta$, then eliminate the choices associated with triples $(i,n-2,n-1)$, since $X_{i,n-2,n-1}$ is not a generator. On the other hand, allow generators of type $X_{i,j,n-2}$ with $a_{n-2} \not\in\Theta$ and $X_{i,j,n-1}$ with $a_{n-1} \not\in\Theta$, for $j\neq n-2$. Then,
\begin{align*}
|A_{3}| &= \binom{n-1-|\Theta|}{3} - \left\{
\begin{array}{cl}
n-3-|\Theta|, & \mbox{if } a_{n-2},a_{n-1} \not\in \Theta; \\ 
0, & \mbox{otherwise}. 
\end{array} 
\right.\\
&= \binom{n-1-|\Theta|}{3} - (1-\mathds{1}_{n-2}^{\Theta})(1-\mathds{1}_{n-1}^{\Theta})(n-3-|\Theta|).
\end{align*}

\item For $|A_{4}\cup A_{5}|$: Count both sets at once. Each connected component of $\Theta$ gives rise to a pair of adjacent roots $a_i$ and $a_{i+1}$ such that $a_i\in \Theta$ and $a_{i+1} \not\in\Theta$, except: 
\begin{enumerate}
\item if there is a connected component which contains in $a_{n-2}$ but not $a_{n-1}$ since it does not exist generator indexed by $(i, n-2,n-1)$;
\item if there is a connected component which contains $a_{n-1}$.
\end{enumerate}
Now, it remains a choice of a root among the remaining $n-1-(|\Theta|-1)=n-2-|\Theta|$ roots to comprise a triple in $A_4 \cup A_5$. Hence,
\begin{align*}
|A_{4}\cup A_{5}| &= \left\{
\begin{array}{cl}
r(n-2-|\Theta|), & \mbox{if } a_{n-2},a_{n-1} \not\in\Theta; \\ 
(r-1)(n-2-|\Theta|), & \mbox{otherwise.}
\end{array} 
\right.\\
&=(r-1+(1-\mathds{1}_{n-2}^{\Theta})(1-\mathds{1}_{n-1}^{\Theta}))(n-2-|\Theta|)\\
&=(r-\mathds{1}_{n-2}^{\Theta}-\mathds{1}_{n-1}^{\Theta}+\mathds{1}_{n-2}^{\Theta}\mathds{1}_{n-1}^{\Theta})(n-2-|\Theta|).
\end{align*}

\item For $|A_{6}|$: Each connected component with two or more roots ending in $a_{i+1}$ gives rise to a sequence of roots $a_i,a_{i+1},a_{i+2}$ such that $a_i,a_{i+1} \in \Theta$ and $a_{i+2}\not\in\Theta$, except the connected component which contains both $a_{n-2}$ and $a_{n-1}$. Hence,
\begin{equation*}
|A_6| = r - r_{1} - \mathds{1}_{n-2}\mathds{1}_{n-1}.
\end{equation*}
\end{itemize}

Therefore, the torsion is easily obtained by the sum $T_{3} = \sum_{k=1}^{6}|A_{k}|-r_{0}$ since $X_{1,1,2}$ generates the free part when $\Theta=\emptyset$. A Sage code to check it is in Appendix \ref{sec:apendice}.
\end{proof}

\begin{theorem}\label{thm:H4}
For $n\geqslant 4$, the dimension of torsion for 4-homology is
\begin{equation*}
T_{4}^{} = \binom{n-|\Theta|+1}{4} + r\binom{n-|\Theta|}{2} + \binom{r}{2}  - (n - |\Theta|-1) (r_{1}+1) - r_{2}.
\end{equation*}

All generators will be presented along the proof.
\end{theorem}
\begin{proof}
Suppose that $\w{i,j,k,l}\in\weyl^{\Theta}$. The boundary map of every $4$-cell in the flag manifold $\mathbb{F}_{\Theta}$ is given by
{\small
\begin{align}
&\textcolor{red}{\partial {\w{i,i,i,i}} = -2(\mathds{1}^{\notInTheta}_{i+2}){\w{i,i+2,i+2}} - 2{\w{i,i,i}}, \  i\in[n-4];^\bullet}\label{eq:b4-1}\\
&\textcolor{blue}{\partial {\w{i,i,i,i+1}} = 0, \  i\in [n-3];^\star}\label{eq:b4-2}\\
&\textcolor{brown}{\partial {\w{i,i,i,i+2}} = -2 {\w{i,i,i+2}}, \  i \in [n-3];^\diamond}\label{eq:b4-3}\\
&\textcolor{red}{\partial {\w{i,i,i,i+3}} = -2 {\w{i,i+2,i+3}} - 2{\w{i,i,i+3}}, \  i \in [n-4];^\bullet}\label{eq:b4-4}\\
&\partial {\w{i,i,i,k}} = -2 (\mathds{1}^{\notInTheta}_{i+2}) {\w{i,i+2,k}},\  i \in [n-5], k\in [i+4,n-1];\label{eq:b4-5}\\
&\textcolor{orange}{\partial {\w{i,i,i+1,i+1}} = 2 (\mathds{1}^{\notInTheta}_{i+2}){\w{i+1,i+1,i+2}}-2 (\mathds{1}^{\notInTheta}_{i}){\w{i,i,i+1}},\  i \in [n-3];^*}\label{eq:b4-6}\\
&\textcolor{brown}{\partial {\w{i,i,i+1,i+2}} = -2{\w{i,i,i+2}},\  i\in [n-3];^\diamond}\label{eq:b4-7}\\
&\textcolor{blue}{\partial {\w{i,i,i+1,j}}= 0,\  i\in [n-4], j \in [i+3,n-1];^\star}\label{eq:b4-8}\\
&\textcolor{brown}{\partial {\w{i,i,i+2,i+2}}= -2{\w{i,i,i+2}},\  i\in [n-4];^\diamond}\label{eq:b4-9}\\
&\textcolor{red}{\partial {\w{i,i,j,j+1}}= -2{\w{i,j,j+1}}-2{\w{i,i,j+1}},\  i\in [n-4], j\in [i+2,n-2];^\bullet}\label{eq:b4-10}\\
&\textcolor{blue}{\partial {\w{i,i,i+2,j}}=0,\  i\in [n-5], j\in [i+4,n-1];^\star}\label{eq:b4-11}\\
&\partial {\w{i,i,j,k}}=-2{\w{i,j,k}},\  i\in [n-6], j\in [i+3,n-3], k\in[j+2,n-1];\label{eq:b4-12}\\
&\textcolor{red}{\partial {\w{i,i,j,j}}=-2{\w{i,j,j}}-2{\w{i,i,j}},\  i\in [n-5],j\in[i+3,n-2];^\bullet}\label{eq:b4-13}\\
&\textcolor{red}{\partial {\w{i,i+1,i+1,i+1}}= -2{\w{i+1,i+1,i+1}} + 2 (\mathds{1}^{\notInTheta}_{i+3}) {\w{i,i+1,i+3}},\  i\in[n-4];^\bullet}\label{eq:b4-14}\\
&\textcolor{blue}{\partial {\w{i,i+1,i+1,i+2}}= 0,\  i\in[n-3];^\star}\label{eq:b4-15}\\
&\textcolor{brown}{\partial {\w{i,i+1,i+1,i+3}}=-2{\w{i+1,i+1,i+3}},\  i\in[n-4];^\diamond}\label{eq:b4-16}\\
&\textcolor{red}{\partial {\w{i,i+1,i+1,j}}=-2{\w{i+1,i+1,j}}+2{\w{i,i+1,j}},\  i\in[n-5], j\in[i+4,n-1];^\bullet}\label{eq:b4-17}\\
&\textcolor{red}{\partial {\w{i,i+1,i+2,i+2}}=2(\mathds{1}^{\notInTheta}_{i}){\w{i,i+2,i+2}}-2{\w{i,i+1,i+2}},\  i\in [n-4];^\bullet}\label{eq:b4-18}\\
&\textcolor{red}{\partial {\w{i,i+1,i+2,i+3}}=-2{\w{i+1,i+2,i+3}}-2(\mathds{1}^{\notInTheta}_{i+1}){\w{i,i+1,i+3}},\  i\in[n-4];^\bullet}\label{eq:b4-19}\\
&\textcolor{red}{\partial {\w{i,i+1,j,j+1}}=-2{\w{i+1,j,j+1}}-2{\w{i,i+1,j+1}},\  i\in[n-5], j\in[i+3,n-2];^\bullet}\label{eq:b4-20}\\
&\textcolor{red}{\partial {\w{i,j,j+1,j+1}}=2{\w{i,j+1,j+1}}-2{\w{i,j,j+1}},\  i\in [n-5], j\in [i+2,n-3];^\bullet}\label{eq:b4-21}\\
&\textcolor{red}{\partial {\w{i,i+1,j,j}}=-2{\w{i+1,j,j}}-2{\w{i,i+1,j}},\  i\in[n-5],j\in[i+3,n-2];^\bullet}\label{eq:b4-22}\\
&\partial {\w{i,i+1,j,k}}=-2{\w{i+1,j,k}},\  i\in [n-6], j\in[i+3,n-3], k\in [j+2,n-1];\label{eq:b4-23}\\
&\partial {\w{i,j,j,j}}=2(\mathds{1}^{\notInTheta}_{j+2}){\w{i,j,j+2}},\  i\in [n-5],j\in [i+2,n-3];\label{eq:b4-24}\\
&\textcolor{blue}{\partial {\w{i,j,j,j+1}}=0,\  i\in[n-4], j\in [i+2,n-2];^\star}\label{eq:b4-25}\\
&\textcolor{blue}{\partial {\w{i,j,j,j+2}}=0,\  i\in[n-5], j\in [i+2,n-3];^\star}\label{eq:b4-26}\\
&\partial {\w{i,j,j,k}}=2{\w{i,j,k}},\  i\in[n-6], j\in [i+2,n-4], k\in[j+3,n-1];\label{eq:b4-27}\\
&\partial {\w{i,j,k,k}}=-2{\w{i,j,k}},\  i\in[n-6], j \in [i+2,n-4], k\in[j+2,n-2];\label{eq:b4-28}\\
&\partial {\w{i,i+1,i+2,k}}=2 (\mathds{1}^{\notInTheta}_{i}) {\w{i,i+2,k}},\  i\in [n-5], k\in[i+4,n-1];\label{eq:b4-29}\\
&\partial {\w{i,j,j+1,k}}=2 {\w{i,j+1,k}},\  i\in [n-6], j \in [i+2,n-4], k\in[j+3,n-1];\label{eq:b4-30}\\
&\partial {\w{i,j,j+1,j+2}}=-2(\mathds{1}^{\notInTheta}_{j}){\w{i,j,j+2}},\  i\in[n-5],j\in[i+2,n-3];\label{eq:b4-31}\\
&\partial {\w{i,j,k,k+1}}=-2{\w{i,j,k+1}},\  i\in[n-6],j\in[i+2,n-4], k\in [j+2,n-2];\label{eq:b4-32}\\
&\textcolor{blue}{\partial {\w{i,j,k,l}}=0,\  i\in [n-7], j\in[i+2,n-5], k\in [j+2,n-3], l\in [k+2,n-1].^\star}\label{eq:b4-33}
\end{align}
}

The generators of $\ker(\partial)$ will be computed by combining the boundary map in Equations \eqref{eq:b4-1} to \eqref{eq:b4-33}. The aim is to get generators of the form $X_{i,j,k,l}$, where $i\leqslant j\leqslant k<l$, in a way that the cell $\w{i,j,k,l}$ is part of the formulas. There are few exceptions to this construction.

The formulas above suggest that there are five kinds of boundary maps (marked with different symbols):

Equation \eqref{eq:b4-6} marked with $*$ provides the free part. Indeed, for $i\in [n-3]$, suppose that $\w{i,i,i+1,i+1}\in \weyl^{\Theta}$, i.e., $a_{i+1}\not\in\Theta$. According to \Cref{tbl:homology_free}, this case will provide all generators $Z_p$, $p\in[r-1]$, for the free part.

The equations marked with $\star$ provide the following generators:
\begin{align*}
&X_{i,i,i,i+1} = {\w{i,i,i,i+1}}, \   i\in [n-3];\\
&X_{i,i,i+i,j} = {\w{i,i,i+1,j}}, \   i\in [n-4], j \in [i+3,n-1];\\
&X_{i,i,i+2,j} = {\w{i,i,i+2,j}}, \   i\in [n-5], j\in [i+4,n-1];\\
&X_{i,i+1,i+1,i+2} = {\w{i,i+1,i+1,i+2}}, \   i\in[n-3];\\
&X_{i,j,j,j+1} = {\w{i,j,j,j+1}}, \   i\in[n-4], j\in [i+2,n-2];\\
&X_{i,j,j,j+2} = {\w{i,j,j,j+2}}, \   i\in[n-5], j\in [i+2,n-3];\\
&X_{i,j,k,l} =  {\w{i,j,k,l}}, \   i\in [n-7], j\in[i+2,n-5],k\in [j+2,n-3],l\in [k+2,n-1].
\end{align*}

For the equations marked with $\diamond$, rewrite Equation \eqref{eq:b4-16} as:
\begin{equation}
\textcolor{brown}{\partial {\w{i-1,i,i,i+2}}=-2{\w{i,i,i+2}},\  i\in[2,n-3].}\tag{\ref*{eq:b4-16}$'$}\label{eq:b4-16p}
\end{equation}

\Cref{tbl:h4_1} summarizes the choices of $i$ such that the boundaries can be combined with each other.
\begin{table}[ht]
\centering
    \caption{Comparing equations marked with $\diamond$}\label{tbl:h4_1}
\begin{tabular}{cc}
\toprule
$i = 1$ & \eqref{eq:b4-3} \eqref{eq:b4-7} \eqref{eq:b4-9}\\ 
$i\in [2,n-4]$  & \eqref{eq:b4-3} \eqref{eq:b4-7} \eqref{eq:b4-9} \eqref{eq:b4-16p} \\  
$i=n-3$ & \eqref{eq:b4-3} \eqref{eq:b4-7} \eqref{eq:b4-16p} \\ 
\bottomrule
\end{tabular} 
\end{table}

These equations provide the following generators:
{\small
\begin{align*}
\mbox{\eqref{eq:b4-3}}-\mbox{\eqref{eq:b4-9}: } & X_{i,i,i,i+2} = {\w{i,i,i,i+2}} - {\w{i,i,i+2,i+2}},\   i\in [n-4];\\
\mbox{\eqref{eq:b4-7}}-\mbox{\eqref{eq:b4-9}: } & X_{i,i,i+1,i+2} = {\w{i,i,i+1,i+2}} - {\w{i,i,i+2,i+2}},\  i\in [n-4];\\
\mbox{\eqref{eq:b4-16p}}-\mbox{\eqref{eq:b4-9}: }& X_{i-1,i,i,i+2} = {\w{i-1,i,i,i+2}} - {\w{i,i,i+2,i+2}},\  i\in [2,n-4]\\
\mbox{or } & X_{i,i+1,i+1,i+3} = {\w{i,i+1,i+1,i+3}} - {\w{i+1,i+1,i+3,i+3}},\  i\in [n-5];\\
\mbox{\eqref{eq:b4-16p}}-\mbox{\eqref{eq:b4-3}: } &X_{n-4,n-3,n-3,n-1} = {\w{n-4,n-3,n-3,n-1}} - {\w{n-3,n-3,n-3,n-1}};\\
\mbox{\eqref{eq:b4-7}}-\mbox{\eqref{eq:b4-3}: } &X_{n-3,n-3,n-2,n-1} = {\w{n-3,n-3,n-2,n-1}} - {\w{n-3,n-3,n-3,n-1}}.
\end{align*}
}

Notice that $X_{n-3,n-3,n-3,n-1}$ is not a generator.

For equations marked with $\bullet$, rewrite some equations as:
\begin{align}
&\textcolor{red}{\partial {\w{i,i,j-1,j}}= -2{\w{i,j-1,j}}-2{\w{i,i,j}},\  i\in [n-4], j\in [i+3,n-1];}\tag{\ref*{eq:b4-10}$'$}\label{eq:b4-10p}\\
&\textcolor{red}{\partial {\w{i-1,i,i,i}}= -2{\w{i,i,i}} + 2 (\mathds{1}^{\notInTheta}_{i+2}) {\w{i-1,i,i+2}},\  i\in[2,n-3];}\tag{\ref*{eq:b4-14}$'$}\label{eq:b4-14p}\\
&\textcolor{red}{\partial {\w{i-1,i,i,j}}=-2{\w{i,i,j}}+2{\w{i-1,i,j}},\  i\in[2,n-4], j\in[i+3,n-1];}\tag{\ref*{eq:b4-17}$'$}\label{eq:b4-17p}\\
&\textcolor{red}{\partial {\w{i-1,i,i+1,i+2}}=-2{\w{i,i+1,i+2}}-2(\mathds{1}^{\notInTheta}_{i}){\w{i-1,i,i+2}},\  i\in[2,n-3];}\tag{\ref*{eq:b4-19}$'$}\label{eq:b4-19p}\\
&\textcolor{red}{\partial {\w{i-1,i,j-1,j}}=-2{\w{i,j-1,j}}-2{\w{i-1,i,j}},\  i\in[2,n-4], j\in[i+3,n-1];}\tag{\ref*{eq:b4-20}$'$}\label{eq:b4-20p}\\
&\textcolor{red}{\partial {\w{i,j-1,j,j}}=2{\w{i,j,j}}-2{\w{i,j-1,j}},\  i\in [n-5], j\in [i+3,n-2];}\tag{\ref*{eq:b4-21}$'$}\label{eq:b4-21p}\\
&\textcolor{red}{\partial {\w{i-1,i,j,j}}=-2{\w{i,j,j}}-2{\w{i-1,i,j}},\  i\in[2,n-4],j\in[i+2,n-2];}\tag{\ref*{eq:b4-22}$'$}\label{eq:b4-22p}
\end{align}

\Cref{tbl:h4_3,tbl:h4_4} summarize the choices of $i$ and $j$ such that the boundaries can be combined with each other, considering that, in some cases, it would be necessary to add up three cells.

\begin{table}[ht]
\centering
    \caption{Equations containing $\w{i,i,i}$}\label{tbl:h4_3}
    \begin{tabular}{cc}
    \toprule
    $i = 1$ & \eqref{eq:b4-1} \\ 

    $i\in [2,n-4]$  & \eqref{eq:b4-1} \eqref{eq:b4-14p} \\ 

    $i=n-3$ & \eqref{eq:b4-14p} \\ 
    \bottomrule 
    \end{tabular} 
\end{table}

\begin{table}[ht]
\centering
    \caption{Comparing equations marked with $\bullet$ containing the given Schubert cell}\label{tbl:h4_4}
    \begin{tabular}{cccccc}
        \toprule
        $i=1$ & $j=3$ & $j=4$ & $j\in [5,n-2]$ & $j = n-1$ \\
        \midrule
        $\w{i,i,j}$ &  & \eqref{eq:b4-4} \eqref{eq:b4-10p} \eqref{eq:b4-13} & \eqref{eq:b4-10p} \eqref{eq:b4-13} & \eqref{eq:b4-10p}  \\ 
        
        $\w{i,j,j}$ & \eqref{eq:b4-1} \eqref{eq:b4-18}  & \eqref{eq:b4-13} \eqref{eq:b4-21p} & \eqref{eq:b4-13} \eqref{eq:b4-21p} & \\ 
    
        $\w{i,j-1,j}$ & \eqref{eq:b4-18} & \eqref{eq:b4-4} \eqref{eq:b4-10p} \eqref{eq:b4-21p}  & \eqref{eq:b4-10p} \eqref{eq:b4-21p} & \eqref{eq:b4-10p} \\ 
        \bottomrule
    
        $i\in [2,n-6]$ & $j=i+2$ & $j=i+3$ & $j\in [i+4,n-2]$ & $j = n-1$ \\
        \midrule
        $\w{i,i,j}$ &  & \eqref{eq:b4-4} \eqref{eq:b4-10p} \eqref{eq:b4-13} \eqref{eq:b4-17p} & \eqref{eq:b4-10p} \eqref{eq:b4-13} \eqref{eq:b4-17p} & \eqref{eq:b4-10p} \eqref{eq:b4-17p} \\ 
      
        $\w{i-1,i,j}$ & \eqref{eq:b4-14p} \eqref{eq:b4-19p} \eqref{eq:b4-22p} & \eqref{eq:b4-17p} \eqref{eq:b4-20p} \eqref{eq:b4-22p} & \eqref{eq:b4-17p} \eqref{eq:b4-20p} \eqref{eq:b4-22p} & \eqref{eq:b4-17p} \eqref{eq:b4-20p}  \\
        
        $\w{i,j,j}$ & \eqref{eq:b4-1} \eqref{eq:b4-18} \eqref{eq:b4-22p} & \eqref{eq:b4-13} \eqref{eq:b4-21p} \eqref{eq:b4-22p} & \eqref{eq:b4-13} \eqref{eq:b4-21p} \eqref{eq:b4-22p}  & \\ 
       
        $\w{i,j-1,j}$ & \eqref{eq:b4-18} \eqref{eq:b4-19p} & \eqref{eq:b4-4} \eqref{eq:b4-10p} \eqref{eq:b4-20p} \eqref{eq:b4-21p}  & \eqref{eq:b4-10p} \eqref{eq:b4-20p} \eqref{eq:b4-21p} & \eqref{eq:b4-10p} \eqref{eq:b4-20p}  \\ 
        \bottomrule
        
        $i = n-5$ & $j=n-3$ & $j=n-2$ & $j = n-1$ \\
        \midrule
        $\w{i,i,j}$ &  & \eqref{eq:b4-4} \eqref{eq:b4-10p} \eqref{eq:b4-13} \eqref{eq:b4-17p} & \eqref{eq:b4-10p} \eqref{eq:b4-17p} \\ 
        
        $\w{i-1,i,j}$ & \eqref{eq:b4-14p} \eqref{eq:b4-19p} \eqref{eq:b4-22p} & \eqref{eq:b4-17p} \eqref{eq:b4-20p} \eqref{eq:b4-22p} & \eqref{eq:b4-17p} \eqref{eq:b4-20p}  \\
        
        $\w{i,j,j}$ & \eqref{eq:b4-1} \eqref{eq:b4-18} \eqref{eq:b4-22p} & \eqref{eq:b4-13} \eqref{eq:b4-21p} \eqref{eq:b4-22p} & \\ 
        
        $\w{i,j-1,j}$ & \eqref{eq:b4-18} \eqref{eq:b4-19p} & \eqref{eq:b4-4} \eqref{eq:b4-10p} \eqref{eq:b4-20p} \eqref{eq:b4-21p}  & \eqref{eq:b4-10p} \eqref{eq:b4-20p}  \\ 

        \bottomrule
        
        $i = n-4$ & $j=n-2$ & $j=n-1$ \\
        \midrule
        
        $\w{i,i,j}$ & & \eqref{eq:b4-4} \eqref{eq:b4-10p} \eqref{eq:b4-17p} \\ 
        
        $\w{i-1,i,j}$ & \eqref{eq:b4-14p} \eqref{eq:b4-19p} \eqref{eq:b4-22p} & \eqref{eq:b4-17p} \eqref{eq:b4-20p} \eqref{eq:b4-22p} \\ 
        
        $\w{i,j,j}$ & \eqref{eq:b4-1} \eqref{eq:b4-18} \eqref{eq:b4-22p} &  \eqref{eq:b4-22p} \\
        
        $\w{i,j-1,j}$ & \eqref{eq:b4-18} \eqref{eq:b4-19p} & \eqref{eq:b4-4} \eqref{eq:b4-10p} \eqref{eq:b4-20p} \\
        \bottomrule
        
        $i=n-3$ & $j= n-1$ \\
        \midrule
        $\w{i-1,i,j}$ & \eqref{eq:b4-14p} \eqref{eq:b4-19p} \\ 
        
        $\w{i,j-1,j}$ & \eqref{eq:b4-19p} \\ 
        \bottomrule
    \end{tabular}
\end{table}

These equations provide the following generators:
{\footnotesize
\begin{align*}
\mbox{\eqref{eq:b4-1}}-\mbox{\eqref{eq:b4-14p}}-\mbox{\eqref{eq:b4-22p}: } &X_{i-1,i,i,n-1} = {\w{i,i,i,i}} - {\w{i-1,i,i,i}} - (\mathds{1}^{\notInTheta}_{i+2}){\w{i-1,i,i+2,i+2}},\  i\in [2,n-4];\\
\mbox{\eqref{eq:b4-4}}-\mbox{\eqref{eq:b4-10p}: } &X_{i,i,i,i+3} ={\w{i,i,i,i+3}} - {\w{i,i,i+2,i+3}}, \  i\in [n-4];\\
\mbox{\eqref{eq:b4-10p}}-\mbox{\eqref{eq:b4-13}}-\mbox{\eqref{eq:b4-21p}: } &X_{i,i,j-1,j} = {\w{i,i,j-1,j}}-{\w{i,i,j,j}}-{\w{i,j-1,j,j}},\  i\in[n-5],j\in[i+3,n-2];\\
\mbox{\eqref{eq:b4-20p}}+\mbox{\eqref{eq:b4-17p}}-\mbox{\eqref{eq:b4-10p}: } &X_{i-1,i,j-1,j} = {\w{i-1,i,j-1,j}} + {\w{i-1,i,i,j}}- {\w{i,i,j-1,j}},\  i\in [2, n-4],j\in [i+3,n-1];\\
\mbox{\eqref{eq:b4-17p}} - \mbox{\eqref{eq:b4-13}}+\mbox{\eqref{eq:b4-22p}: } &X_{i-1,i,i,j} = {\w{i-1,i,i,j}} - {\w{i,i,j,j}} + {\w{i-1,i,j,j}},\  i\in [2,n-5],j\in [i+3,n-2];\\
\mbox{\eqref{eq:b4-19p}}-\mbox{\eqref{eq:b4-18}}-\mbox{\eqref{eq:b4-22p}: } &X_{i-1,i,i+1,i+2} = {\w{i-1,i,i+1,i+2}} - {\w{i,i+1,i+2,i+2}} - (\mathds{1}^{\notInTheta}_{i}){\w{i-1,i,i+2,i+2}},\  i\in [2,n-4];\\
\mbox{\eqref{eq:b4-20p}}-\mbox{\eqref{eq:b4-21p}}-\mbox{\eqref{eq:b4-22p}: } &\tilde X_{i,j} = {\w{i-1,i,j-1,j}} - {\w{i,j-1,j,j}} - {\w{i-1,i,j,j}},\  i \in [2,n-5],j\in [i+3,n-2].
\end{align*}
}

Observe that $\tilde X_{i,j}$ is not a generator since $\tilde X_{i,j} = X_{i,i,j-1,j} + X_{i-1,i,j-1,j} - X_{i-1,i,i,j}$ when $i\in [2,n-5]$, $j\in [i+3,n-2]$.

For the equations with no marker, rewrite Equations \eqref{eq:b4-23}, \eqref{eq:b4-30}, and \eqref{eq:b4-32} as:
\begin{align}
&\partial {\w{i-1,i,j,k}}=-2{\w{i,j,k}},\  i\in [2,n-5], j\in[i+2,n-3], k\in [j+2,n-1];\tag{\ref*{eq:b4-23}$'$}\label{eq:b4-23p}\\
&\partial {\w{i,j,j+1,k}}=2 {\w{i,j+1,k}},\  i\in [n-6], j \in [i+2,n-4], k\in[j+3,n-1];\tag{\ref*{eq:b4-30}$'$}\label{eq:b4-30p}\\
&\partial {\w{i,j,k-1,k}}=-2{\w{i,j,k}},\  i\in[n-6],j\in[i+2,n-4], k\in [j+3,n-1].\tag{\ref*{eq:b4-32}$'$}\label{eq:b4-32p}
\end{align}

\Cref{tbl:h4_2} summarize the choices of $i$, $j$, and $k$ such that the boundaries can be combined with each other.

\begin{table}[ht]
\centering
    \caption{Comparing equations without marker}\label{tbl:h4_2}
    \begin{tabular}[c]{cccc}
        \toprule
         $i=1$ & $k=j+2$ & $k\in [j+3,n-2]$ & $k=n-1$ \\
        \midrule
        $j=i+2$ & \eqref{eq:b4-5} \eqref{eq:b4-24} \eqref{eq:b4-28} \eqref{eq:b4-29} \eqref{eq:b4-31} & \eqref{eq:b4-5} \eqref{eq:b4-27} \eqref{eq:b4-28} \eqref{eq:b4-29} \eqref{eq:b4-32p} & \eqref{eq:b4-5} \eqref{eq:b4-27} \eqref{eq:b4-29} \eqref{eq:b4-32p} \\
        
        $j\in [i+3, n-4]$ & \eqref{eq:b4-12} \eqref{eq:b4-24} \eqref{eq:b4-28} \eqref{eq:b4-30p} \eqref{eq:b4-31} & \eqref{eq:b4-12} \eqref{eq:b4-27} \eqref{eq:b4-28} \eqref{eq:b4-30p} \eqref{eq:b4-32p} & \eqref{eq:b4-12} \eqref{eq:b4-27} \eqref{eq:b4-30p} \eqref{eq:b4-32p} \\ 
        
        $j=n-3$ & \eqref{eq:b4-12} \eqref{eq:b4-24} \eqref{eq:b4-30p} \eqref{eq:b4-31} & \eqref{eq:b4-12} \eqref{eq:b4-30p} & \eqref{eq:b4-12} \eqref{eq:b4-30p} \\ 
        \bottomrule

        $i\in [2,n-6]$ & $k=j+2$ & $k\in [j+3,n-2]$ & $k=n-1$ \\
        \midrule
        $j=i+2$ & \eqref{eq:b4-5} \eqref{eq:b4-23p} \eqref{eq:b4-24} \eqref{eq:b4-28} \eqref{eq:b4-29} \eqref{eq:b4-31} & \eqref{eq:b4-5} \eqref{eq:b4-23p} \eqref{eq:b4-27} \eqref{eq:b4-28} \eqref{eq:b4-29} \eqref{eq:b4-32p} & \eqref{eq:b4-5} \eqref{eq:b4-23p} \eqref{eq:b4-27} \eqref{eq:b4-29} \eqref{eq:b4-32p} \\ 
        
        $j\in [i+3, n-4]$ & \eqref{eq:b4-12} \eqref{eq:b4-23p} \eqref{eq:b4-24} \eqref{eq:b4-28} \eqref{eq:b4-30p} \eqref{eq:b4-31} & \eqref{eq:b4-12} \eqref{eq:b4-23p} \eqref{eq:b4-27} \eqref{eq:b4-28} \eqref{eq:b4-30p} \eqref{eq:b4-32p} & \eqref{eq:b4-12} \eqref{eq:b4-23p} \eqref{eq:b4-27} \eqref{eq:b4-30p} \eqref{eq:b4-32p}\\ 
        
        $j=n-3$ & \eqref{eq:b4-12} \eqref{eq:b4-23p} \eqref{eq:b4-24} \eqref{eq:b4-30p} \eqref{eq:b4-31} & \eqref{eq:b4-12} \eqref{eq:b4-23p} \eqref{eq:b4-30p} & \eqref{eq:b4-12} \eqref{eq:b4-23p} \eqref{eq:b4-30p} \\ 
        \bottomrule
        
        $i=n-5$ & $k=j+2$ & $k\in [j+3,n-2]$ & $k=n-1$ \\
        \midrule
        $j=i+2$ & \eqref{eq:b4-5} \eqref{eq:b4-23p} \eqref{eq:b4-24} \eqref{eq:b4-29} \eqref{eq:b4-31} & \eqref{eq:b4-5} \eqref{eq:b4-23p} \eqref{eq:b4-29} & \eqref{eq:b4-5} \eqref{eq:b4-23p} \eqref{eq:b4-29} \\ 
        
        $j\in [i+3, n-4]$ & \eqref{eq:b4-23p} \eqref{eq:b4-24} \eqref{eq:b4-31} & \eqref{eq:b4-23p} & \eqref{eq:b4-23p} \\ 
        
        $j=n-3$ & \eqref{eq:b4-23p} \eqref{eq:b4-24} \eqref{eq:b4-31} & \eqref{eq:b4-23p} & \eqref{eq:b4-23p}  \\ 
        \bottomrule
    \end{tabular}
\end{table}

These equations provide the following generators:
{\footnotesize
\begin{align*}
\mbox{\eqref{eq:b4-5}}-\mbox{\eqref{eq:b4-28}: } & X_{1,1,1,5} = {\w{1,1,1,5}}-(\mathds{1}^{\notInTheta}_{3}) {\w{1,3,5,5}};\\
\mbox{\eqref{eq:b4-24}}+\mbox{\eqref{eq:b4-28}: } & Y_{1,3} = {\w{1,3,3,3}}+(\mathds{1}^{\notInTheta}_{5}) {\w{1,3,5,5}};\\
\mbox{\eqref{eq:b4-29}}+\mbox{\eqref{eq:b4-28}: } &X_{1,2,3,5} = {\w{1,2,3,5}}+(\mathds{1}^{\notInTheta}_{1}) {\w{1,3,5,5}};\\
\mbox{\eqref{eq:b4-31}}-\mbox{\eqref{eq:b4-28}: } &X_{1,3,4,5} = {\w{1,3,4,5}}-(\mathds{1}^{\notInTheta}_{3}) {\w{1,3,5,5}};\\
\mbox{\eqref{eq:b4-5}}-\mbox{\eqref{eq:b4-32p}: } & X_{1,1,1,k} = {\w{1,1,1,k}}-(\mathds{1}^{\notInTheta}_{3}){\w{1,3,k-1,k}}, \  k\in [6,n-1],(\mbox{take } i = 1 ,j = i+2 = 3);\\
\mbox{\eqref{eq:b4-27}}+\mbox{\eqref{eq:b4-32p}: } &X_{1,3,3,k} = {\w{1,3,3,k}}+{\w{1,3,k-1,k}},\  k\in[6,n-1],(\mbox{take } i = 1, j=3);\\
\mbox{\eqref{eq:b4-32p}}-\mbox{\eqref{eq:b4-28}: } &X_{1,3,k-1,k} = {\w{1,3,k-1,k}}-{\w{1,3,k,k}},\  k\in[6,n-2],(\mbox{take } i = 1);\\
\mbox{\eqref{eq:b4-29}}+\mbox{\eqref{eq:b4-32p}: } &X_{1,2,3,k} = {\w{1,2,3,k}}+(\mathds{1}^{\notInTheta}_{1}){\w{1,3,k-1,k}}, \  k\in [6,n-1],(\mbox{take } i = 1 ,j = i+2 = 3);\\
\mbox{\eqref{eq:b4-12}}+\mbox{\eqref{eq:b4-30p}: } & X_{1,1,j,k} = {\w{1,1,j,k}}+{\w{1,j-1,j,k}},\  j\in [4,n-3],k\in[j+2,n-1],(\mbox{take } i = 1);\\
\mbox{\eqref{eq:b4-24}}-\mbox{\eqref{eq:b4-30p}: } & Y_{1,j} = {\w{1,j,j,j}}-(\mathds{1}^{\notInTheta}_{j+2}){\w{1,j-1,j,j+2}},\  j\in [4,n-3],(\mbox{take } i = 1,k = j+2);\\
\mbox{\eqref{eq:b4-27}}-\mbox{\eqref{eq:b4-30p}: } & X_{1,j,j,k} = {\w{1,j,j,k}}-{\w{1,j-1,j,k}}, \  j\in [4,n-4],k\in [j+3,n-1],(\mbox{take } i = 1);\\
\mbox{\eqref{eq:b4-28}}+\mbox{\eqref{eq:b4-30p}: } &X_{1,j-1,j,k} = {\w{1,j-1,j,k}}+{\w{1,j,k,k}},\  j\in [4,n-4],k\in [j+2,n-2],(\mbox{take } i = 1);\\
\mbox{\eqref{eq:b4-31}}+\mbox{\eqref{eq:b4-30p}: } &X_{1,j,j+1,j+2} = {\w{1,j,j+1,j+2}}+(\mathds{1}^{\notInTheta}_{j}){\w{1,j-1,j,j+2}},\  j\in [4,n-3],(\mbox{take } i = 1 ,k = j+2);\\
\mbox{\eqref{eq:b4-30p}}+\mbox{\eqref{eq:b4-32p}: } &X_{1,j,k-1,k} = {\w{1,j,k-1,k}}+{\w{1,j-1,j,k}},\  j\in [4,n-4],k\in[j+3,n-1],(\mbox{take } i = 1);\\
\mbox{\eqref{eq:b4-5}}-\mbox{\eqref{eq:b4-23p}: } &X_{i,i,i,k} = {\w{i,i,i,k}}-(\mathds{1}^{\notInTheta}_{i+2}){\w{i-1,i,i+2,k}},\  i\in [2,n-5],k\in [i+4,n-1],(\mbox{take } j = i+2);\\
\mbox{\eqref{eq:b4-12}}-\mbox{\eqref{eq:b4-23p}: } &X_{i,i,j,k} = {\w{i,i,j,k}}-{\w{i-1,i,j,k}},\  i\in [2,n-6],j\in [i+3,n-3],k\in[j+2,n-1];\\
\mbox{\eqref{eq:b4-24}}+\mbox{\eqref{eq:b4-23p}: } & Y_{i,j} = {\w{i,j,j,j}}+(\mathds{1}^{\notInTheta}_{j+2}){\w{i-1,i,j,j+2}},\  i\in [2,n-5],j\in [i+2,n-3],(\mbox{take } k = j+2);\\
\mbox{\eqref{eq:b4-27}}+\mbox{\eqref{eq:b4-23p}: } &X_{i,j,j,k} = {\w{i,j,j,k}}+{\w{i-1,i,j,k}},\  i\in [2,n-6],j\in [i+2,n-4],k\in[j+3,n-1];\\
\mbox{\eqref{eq:b4-23p}}-\mbox{\eqref{eq:b4-28}: } &X_{i-1,i,j,k} = {\w{i-1,i,j,k}}-{\w{i,j,k,k}},\  i\in [2,n-6],j\in [i+2,n-4],k\in[j+2,n-2];\\
\mbox{\eqref{eq:b4-29}}+\mbox{\eqref{eq:b4-23p}: } &X_{i,i+1,i+2,k} = {\w{i,i+1,i+2,k}}+(\mathds{1}^{\notInTheta}_{i}){\w{i-1,i,i+2,k}},\  i\in [2,n-5],k\in [i+4,n-1],(\mbox{take } j = i+2);\\
\mbox{\eqref{eq:b4-30p}}+\mbox{\eqref{eq:b4-23p}: } &X_{i,j-1,j,k} = {\w{i,j-1,j,k}}+{\w{i-1,i,j,k}},\  i\in [2,n-6],j\in [i+3,n-3],k\in[j+2,n-1];\\
\mbox{\eqref{eq:b4-31}}+\mbox{\eqref{eq:b4-23p}: } &X_{i,j,j+1,j+2} = {\w{i,j,j+1,j+2}}+(\mathds{1}^{\notInTheta}_{j}){\w{i-1,i,j,j+2}},\  i\in [2,n-5],j\in [i+2,n-3],(\mbox{take } k = j+2);\\
\mbox{\eqref{eq:b4-32p}}+\mbox{\eqref{eq:b4-23p}: } &X_{i,j,k-1,k} = {\w{i,j,k-1,k}}+{\w{i-1,i,j,k}},\  i\in [2,n-6],j\in [i+2,n-4],k\in[j+3,n-1].
\end{align*}
}

In order to use a generator of the form $Y_{i,j}$, where $i\leqslant j+2\leqslant n-1$, rewrite them by
\begin{align*}
& X_{1,3,n-2,n-1} = Y_{1,3};\\
& X_{1,j-1,j,n-1} = Y_{1,j}, \  j\in[4,n-3];\\
& X_{i-1,i,j,n-1} = Y_{i,j},\  i\in[2,n-5], j\in[i+2,n-3].
\end{align*}

Notice in this case that either $n-1\in\Theta$ or $n-1\not\in\Theta$.

Except for $Z_i$, all other kernel generators are torsion generators. The torsion generators will be counted by splitting up according to the following disjoint sets:
{\small
\begin{align*}
    A_1 &=\{X_{i,i,i,k} \colon a_i,a_{k} \not\in\Theta\}; & A_9 &=\{X_{i-1,i,j,k} \colon a_{i},a_{j},a_{k} \not\in\Theta; a_{i-1}\in\Theta\};\\
    A_2 &=\{X_{i,i,j,k} \colon a_{i},a_{j},a_{k} \not\in\Theta\}; & A_{10} &=\{X_{i,j-1,j,k} \colon a_{i},a_{j},a_{k} \not\in\Theta; a_{j-1}\in\Theta\};\\
    A_3 &=\{X_{i,j,j,k} \colon a_{i},a_{j},a_{k} \not\in\Theta\}; & A_{11} &=\{X_{i,j,k-1,k} \colon a_{i},a_{j},a_{k} \not\in\Theta; a_{k-1}\in\Theta\};\\
    A_4 &=\{X_{i,j,k,l} \colon a_{i},a_{j},a_{k},a_{l} \not\in\Theta\}; & A_{12} &=\{X_{i,j,k,n-1} \colon a_{i},a_{j},a_{k} \not\in\Theta; a_{n-1}\in\Theta\};\\
    A_5 &=\{X_{i,i,j-1,j} \colon a_{i},a_{j} \not\in\Theta; a_{j-1}\in\Theta\}; & A_{13} &=\{X_{i-1,i,j-1,j} \colon a_{i},a_{j}\not\in\Theta; a_{i-1},a_{j-1}\in\Theta\};\\
    A_6 &=\{X_{i-1,i,i,j} \colon a_{i},a_{j} \not\in\Theta; a_{i-1}\in\Theta\}; & A_{14} &=\{X_{i-1,i,i+1,j} \colon a_{i+1},a_{j}\not\in\Theta; a_{i-1},a_{i}\in\Theta\};\\
    A_7 &=\{X_{i-1,i,i,n-1} \colon a_{i-1},a_{i} \not\in\Theta; a_{n-1}\in\Theta\}; & A_{15} &=\{X_{i,j-1,j,j+1} \colon a_{i},a_{j+1}\not\in\Theta; a_{j-1},a_{j}\in\Theta\};\\
    A_{8} &=\{X_{i-1,i,i,n-1} \colon a_{i} \not\in\Theta; a_{i-1},a_{n-1}\in\Theta\}; & A_{16} &=\{X_{i-1,i,i+1,i+2} \colon a_{i+2}\not\in\Theta; a_{i-1},a_{i},a_{i+1}\in\Theta\}.
\end{align*}
}

It follows that:
\begin{itemize}
\item For $|A_{1}|$: There is no generator indexed by $(n-2,n-2,n-2,n-1)$ and $(n-3,n-3,n-3,n-1)$. Then,
\begin{align*}
|A_{1}| & = \binom{n-1-|\Theta|}{2} -
\left\{
\begin{array}{cl}
0, & \mbox{if } a_{n-1}\in\Theta; \\ 
0, & \mbox{if } a_{n-3},a_{n-2}\in\Theta \mbox{ and } a_{n-1}\not\in\Theta; \\
1, & \mbox{if } a_{n-3}\in\Theta \mbox{ and } a_{n-2},a_{n-1}\not\in\Theta; \\
1, & \mbox{if } a_{n-2}\in\Theta \mbox{ and } a_{n-3},a_{n-1}\not\in\Theta; \\
2, & \mbox{if } a_{n-3}, a_{n-2},a_{n-1}\not\in\Theta.
\end{array} 
\right. \\
& =\binom{n-1-|\Theta|}{2} - (2-\one^{\Theta}_{n-2}-\one^{\Theta}_{n-3})(1-\one^{\Theta}_{n-1}).
\end{align*}

\item For $|A_{2}|$: There is no generator indexed by $(i,i,n-2,n-1)$, for $i\in[n-4]$. Then,
\begin{align*}
|A_{2}| & = \binom{n-1-|\Theta|}{3} -
\left\{
\begin{array}{cl}
0, & \mbox{if } a_{n-1}\in\Theta \mbox{ or } a_{n-2}\in\Theta; \\
n-4-|\Theta|, & \mbox{if } a_{n-3},a_{n-2},a_{n-1}\not\in\Theta; \\
n-3-|\Theta|, & \mbox{if } a_{n-3}\in\Theta \mbox{ and } a_{n-2},a_{n-1}\not\in\Theta.
\end{array} 
\right. \\
& =\binom{n-1-|\Theta|}{3} - (1-\one^{\Theta}_{n-2})(1-\one^{\Theta}_{n-1}) (n-4 - |\Theta| +\one^{\Theta}_{n-3}).
\end{align*}

\item For $|A_{3}|$: Then,
\begin{align*}
|A_{3}| & = \binom{n-1-|\Theta|}{3}.
\end{align*}

\item For $|A_{4}|$: There is no generator indexed by $(n-4,n-3,n-2,n-1)$. Then,
\begin{align*}
|A_{4}| & = \binom{n-1-|\Theta|}{4} -
(1-\one^{\Theta}_{n-4})(1-\one^{\Theta}_{n-3})(1-\one^{\Theta}_{n-2})(1-\one^{\Theta}_{n-1}).
\end{align*}

\item For $|A_{5}\cup A_{6}|$: There is no generator indexed by $(i,i,n-2,n-1)$, for $i\in[n-4]$. Then,
\begin{align*}
|A_{5}\cup A_{6}| & = 
\left\{
\begin{array}{cl}
(r-1)(n-2-|\Theta|), & \mbox{if } a_{n-1}\in\Theta; \\
(r-1)(n-2-|\Theta|) + 1, & \mbox{if } a_{n-2}\in\Theta \mbox{ and } a_{n-3}, a_{n-1}\not\in\Theta; \\
(r-1)(n-2-|\Theta|), & \mbox{if } a_{n-3}, a_{n-2}\in\Theta \mbox{ and } a_{n-1}\not\in\Theta; \\
r(n-2-|\Theta|), & \mbox{if } a_{n-2}, a_{n-1}\not\in\Theta.
\end{array} 
\right. \\
& = (r-\one^{\Theta}_{n-2}-\one^{\Theta}_{n-1}(1-\one^{\Theta}_{n-2}))(n-2-|\Theta|)+(1-\one^{\Theta}_{n-3})\one^{\Theta}_{n-2}(1-\one^{\Theta}_{n-1}).
\end{align*}

\item For $|A_{7}|$:
Then,
\begin{equation*}
|A_{7}| = \one^{\Theta}_{n-1}\left( \sum_{i=2}^{n-4} (1-\one_{i-1}^{\Theta})(1-\one^{\Theta}_{i})\right).
\end{equation*}

\item For $|A_{8}|$:
Then,
\begin{equation*}
|A_{8}| = \one^{\Theta}_{n-1}\left( \sum_{i=2}^{n-4} \one_{i-1}^{\Theta}(1-\one^{\Theta}_{i})\right)
 = \one_{n-1}^{\Theta}(r-1 - \one_{n-4}^{\Theta}(1-\one_{n-3}^{\Theta}) - \one_{n-3}^{\Theta}(1-\one_{n-2}^{\Theta})),
\end{equation*}
since $\displaystyle\sum_{i=2}^{n-1} \one_{i-1}^{\Theta}(1-\one^{\Theta}_{i}) = r-1$.

\item For $|A_{9}\cup A_{10}\cup A_{11}|$: There is no generator indexed by $(n-4,n-3,n-2,n-1)$. Then,
\begin{align*}
|A_{9}\cup A_{10}\cup A_{11}| = (r-\one^{\Theta}_{n-1}) \binom{n-2-|\Theta|}{2} &- \one^{\Theta}_{n-4} (1-\one^{\Theta}_{n-3}) (1-\one^{\Theta}_{n-2})(1-\one^{\Theta}_{n-1}) \\
& - (1-\one^{\Theta}_{n-4}) \one^{\Theta}_{n-3} (1-\one^{\Theta}_{n-2})(1-\one^{\Theta}_{n-1})\\
&-(1-\one^{\Theta}_{n-4}) (1-\one^{\Theta}_{n-3}) \one^{\Theta}_{n-2}(1-\one^{\Theta}_{n-1}).
\end{align*}

\item For $|A_{12}|$:
Then,
\begin{align*}
|A_{12}|=\one_{n-1}^{\Theta}\left(\binom{n-2-|\Theta|}{2} +\one_{n-2}^{\Theta}(n-2-|\Theta|)-\sum_{i=2}^{n-3} (1-\one^{\Theta}_{i-1})(1-\one^{\Theta}_{i})\right).
\end{align*}

\item For $|A_{13}|$: There is no generator indexed by $(n-4,n-3,n-2,n-1)$. Then,
\begin{align*}
|A_{13}| = \binom{r}{2} - (r-1)\one_{n-1}^{\Theta} - \one^{\Theta}_{n-4} (1-\one^{\Theta}_{n-3}) \one^{\Theta}_{n-2}(1-\one_{n-1}^{\Theta}).
\end{align*}

\item For $|A_{14}\cup A_{15}|$: There is no generator indexed by $(n-4,n-3,n-2,n-1)$. Then,
\begin{align*}
|A_{14}\cup A_{15}| = 
(r - r_{1} - \one_{n-2}^{\Theta}\one_{n-1}^{\Theta})(n-2-|\Theta|) & - \one^{\Theta}_{n-4} \one^{\Theta}_{n-3} (1-\one^{\Theta}_{n-2})(1-\one^{\Theta}_{n-1}) \\
&- (1-\one^{\Theta}_{n-4}) \one^{\Theta}_{n-3} \one^{\Theta}_{n-2}(1-\one^{\Theta}_{n-1}).
\end{align*}

\item For $|A_{16}|$: There is no generator indexed by $(n-4,n-3,n-2,n-1)$. Then,
\begin{align*}
|A_{16}| = 
r-r_{1}-r_{2} - \one^{\Theta}_{n-4} \one^{\Theta}_{n-3} \one^{\Theta}_{n-2}(1-\one_{n-1}^{\Theta}) - \one^{\Theta}_{n-3} \one^{\Theta}_{n-2} \one^{\Theta}_{n-1}.
\end{align*}
\end{itemize}

Therefore, the torsion is obtained by the sum $T_{4} = \sum_{k=1}^{16}|A_{k}|$. A Sage code to verify it is in the Appendix \ref{sec:apendice}.
\end{proof}

\section*{Acknowledgments}

We thank San Martin and Lucas Seco for helpful suggestions and valuable discussions. This research was motivated by computer investigation using the open-source mathematical software Sage \cite{sagemath}.

\bibliographystyle{plain} 
\bibliography{biblio}

\appendix

\section{Sage code}\label{sec:apendice}

We use the following Sage code to verify the final sum in the proof of \Cref{thm:H3}. Denote \texttt{t}=$|\Theta|$, \texttt{r}=$r$, \texttt{r1}=$r_{1}$, \texttt{I[0]}=$\one^{\Theta}_{n-2}$, \texttt{I[1]}=$\one^{\Theta}_{n-1}$.

\begin{lstlisting}[language=Python]
# Sage code to check the torsion coefficient of H_3
var('n,t,r,r1')

def A(x,I):
    if x == 1:
        return (n-3-t) + I[0] + I[1]
    if x == 2:
        return binomial(n-1-t,2)
    if x == 3:
        return binomial(n-1-t,3) - (1-I[0])*(1-I[1])*(n-3-t)
    if x == 4:
        return (r-I[0]-I[1]+I[0]*I[1])*(n-2-t)
    if x == 6:
        return r - r1 - I[0]*I[1]
    else:
        return 0

#Claimed formula
T_3 = binomial(n-t,3) + r*(n-t-1) - r1

# Compute the sum for each possibility of I and compare it with the claimed formula
for i in [0,1]:
    for j in  [0,1]:
        SUM = sum([A(x, [i,j]) for x in range(1,7)])
        
        # Compare the sum with the claimed formula
        print(bool(SUM==T_3))
\end{lstlisting}

We use the following Sage code to verify the final sum in the proof of \Cref{thm:H4}. Denote \texttt{t}=$|\Theta|$, \texttt{r}=$r$, \texttt{r1}=$r_{1}$, \texttt{r2}=$r_{2}$, \texttt{S}=$\sum_{i=2}^{n-4}(1-\one^{\Theta}_{i-1})(1-\one^{\Theta}_{i})$, \texttt{I[0]}=$\one^{\Theta}_{n-4}$, \texttt{I[1]}=$\one^{\Theta}_{n-3}$, \texttt{I[2]}=$\one^{\Theta}_{n-2}$, \texttt{I[3]}=$\one^{\Theta}_{n-1}$.
\begin{lstlisting}[language=Python]
# Sage code to check the torsion coefficient of H_4
var('n,t,r,r1,r2,S')

def A(x,I):
    if x == 1:
        return binomial(n-1-t,2) - (2-I[1]-I[2])*(1-I[3])
    if x == 2:
        return binomial(n-1-t,3) - (1-I[2])*(1-I[3])*(n-4-t+I[1])
    if x == 3:
        return binomial(n-1-t,3)
    if x == 4:
        return binomial(n-1-t,4) - (1-I[0])*(1-I[1])*(1-I[2])*(1-I[3])
    if x == 5:
        return (r-I[2]-I[3]*(1-I[2]))*(n-2-t)+(1-I[1])*I[2]*(1-I[3])
    if x == 7:
        return I[3]*S
    if x == 8:
        return I[3]*(r-1-I[0]*(1-I[1])-I[1]*(1-I[2]))
    if x == 9:
        return (r-I[3])*binomial(n-2-t,2) - I[0]*(1-I[1])*(1-I[2])*(1-I[3]) - (1-I[0])*I[1]*(1-I[2])*(1-I[3]) - (1-I[0])*(1-I[1])*I[2]*(1-I[3])
    if x == 12:
        return I[3]*(binomial(n-2-t,2) + I[2]*(n-2-t) - S - (1-I[0])*(1-I[1]))
    if x == 13:
        return binomial(r,2) - (r-1)*I[3] - I[0]*(1-I[1])*I[2]*(1-I[3])
    if x == 14:
        return (r-r1-I[2]*I[3])*(n-2-t) - I[0]*I[1]*(1-I[2])*(1-I[3]) - (1-I[0])*I[1]*I[2]*(1-I[3])
    if x == 16:
        return r - r1 - r2 - I[0]*I[1]*I[2]*(1-I[3])-I[1]*I[2]*I[3]
    else:
        return 0

#Claimed formula
T_4 = binomial(n-t+1,4) + r*binomial(n-t,2) + binomial(r,2) - (n-t-1)*(r1+1) - r2

# Compute the sum for each possibility of I and compare it with the claimed formula
for i in [0,1]:
    for j in  [0,1]:
        for k in  [0,1]:
            for l in  [0,1]:
                SUM = sum([A(x, [i,j,k,l]) for x in range(1,17)])
                
                # Compare the sum with the claimed formula
                print(bool(SUM==T_4))
\end{lstlisting}

\end{document}